\documentclass[letterpaper,journal]{IEEEtran}
\usepackage[T1]{fontenc}
\usepackage[latin9]{inputenc}
\usepackage{color}
\usepackage{babel}
\usepackage[caption=false,font=normalsize,labelfont=sf,textfont=sf]{subfig}
\usepackage{textcomp}
\usepackage{stfloats}
\usepackage{verbatim}
\usepackage{mathrsfs,bbm}
\usepackage{amsmath}
\usepackage{amsthm}
\usepackage{amssymb,comment}
\usepackage{stackrel}
\usepackage{graphicx}
\usepackage{algorithm} 
\usepackage{algpseudocode} 
\usepackage[hyphens]{url}   

\newcommand{\ubeta}{^{(\beta)}}
\newcommand{\nal}[1]{\begin{align*}#1\end{align*}}
\newcommand{\al}[1]{\begin{align}#1\end{align}}

\usepackage[unicode=true,
 bookmarks=true,bookmarksnumbered=true,bookmarksopen=true,bookmarksopenlevel=1,
 breaklinks=false,pdfborder={0 0 1},backref=false,colorlinks=false]
 {hyperref}
\hypersetup{pdftitle={Your Title},
 pdfauthor={Your Name},
 pdfborderstyle=,pdfpagelayout=OneColumn,pdfnewwindow=true,pdfstartview=XYZ,plainpages=false}

\makeatletter

\usepackage{balance}
\usepackage{color}
\usepackage{ragged2e}

\newcommand{\ust}{^{\star}}
\newcommand{\cF}{\mathcal{F}}

\newcommand{\te}{\theta}
\newcommand{\bR}{\mathbb{R}}
\newcommand{\bE}{\mathbb{E}}
\newcommand{\cU}{\mathcal{U}}

\newcommand{\cQ}{\mathcal{Q}}
\newcommand{\bP}{\mathbb{P}}
\newcommand{\bZ}{\mathbb{Z}}
\newcommand{\cB}{\mathcal{B}}
\newcommand{\cX}{\mathcal{X}}
\newcommand{\cG}{\mathcal{G}}

\newcommand{\bN}{\mathbb{N}}

\newtheorem{proposition}{Proposition}

\newtheorem{definition}{Definition}
\newtheorem{theorem}{Theorem}
\newtheorem{lemma}{Lemma}
\newtheorem{assumption}{Assumption}

\makeatother

\begin{document}
\title{Delay-Optimal Transmission Scheduling Policies for Time-Correlated Fading Channels}
\author{Manali Dutta, Gourav~Saha,~Rahul~Singh,~and~Ness~B.~Shroff\thanks{This material is based upon work supported by the National Science
Foundation under grant number CNS-1955535.\protect \\
 Manali Dutta is with the Department of Computer Science and Automation,
Indian Institute of Science, Bengaluru, Karnataka 560012,
India; Email: manalidutta@iisc.ac.in.\protect\\
Gourav Saha is with the Department of Computer Science and Engineering,
Mahindra University, Hyderabad, Telangana 500043, India; Email: gourav.saha@mahindrauniversity.edu.in.\protect \\
Email: rahulsingh0188@gmail.com.\protect \\
 Ness B. Shroff is with the Department of Electrical
and Computer Engineering, Ohio State University, Columbus, OH 43210,
USA; Email: shroff.11@osu.edu.}}
\maketitle
\begin{abstract}
Millimeter-wave (mmWave) networks have the potential to support high throughput and low-latency requirements of 5G-and-beyond communication standards. But transmissions in this band are highly vulnerable to attenuation and blockages from humans, buildings, and foliage, which increase end-to-end packet delays.
~This work designs dynamic scheduling policies that minimize end-to-end packet delays while keeping packet transmission costs low.
~Specifically, we consider a mmWave network that consists of a transmitter that transmits data packets over an unreliable communication channel modeled as a Gilbert-Elliott channel.
The transmitter operates under an ACK/NACK feedback model and does not observe the channel state unless it attempts a transmission.~The objective is to minimize a weighted average cost consisting of end-to-end packet delays and packet transmission costs.
~We pose this dynamic optimization problem as a partially observable
Markov decision process (POMDP). To the best of our knowledge, this is the first POMDP formulation for mmWave network with partial channel state information that considers delay minimization.~We show that the POMDP admits a solution that has a ``threshold'' structure, i.e., for each queue length, the belief (the conditional probability that the channel is in a good state) is partitioned into intervals, and the transmitter sends $j$ packets when the belief lies in the $j$-th interval.
~We then consider the case when the system parameters such as the packet arrival rate, and the transition probabilities of the channel are not known, and leverage these structural results in order to use the actor-critic algorithm to efficiently search for a policy that is locally optimal.
\end{abstract}

\begin{IEEEkeywords}
Millimeter wave, latency, partially observable Markov decision process,
threshold policy, actor-critic algorithm.\vspace{-1em}
 
\end{IEEEkeywords}

\section{Introduction\label{sec:Introduction}}
By 2030, typical 
global (download) data traffic would require ultra-low latency guarantees of the order of milliseconds, and a throughput that may reach up to several Gbps~\cite{wbb}. 
This demand is driven by the ``5G-and-beyond'' wireless networks that have to support a large number of IoT devices \cite{mehmood2017m2m}.~Indeed, applications such as vehicle-to-vehicle communications, multimedia,
and industrial control networks require low latency communication~\cite{wang2018millimeter,yang2018low}. It
is predicted that the mobile data traffic would
increase by $1000$-folds every decade~\cite{li20145g}.~Since it is critical for Ultra-Reliable Low Latency
Communications (URLLC)~\cite{popovski2017ultra} that the data packets have a low end-to-end latency, such networks impose stringent Quality-of-Service (QoS) guarantees in order to support emerging applications~\cite{popovski2017ultra,wang2018joint,andrews2014will}.~Millimeter-wave
(mmWave) communication networks have the capability
of providing data rates of the order of multi-gigabits per second, and can meet this demand because of the following properties that distinguish it from the sub-6 GHz band networks~\cite{yang2018low}: abundancy of spectral resources, low
costs of component, and the use of highly directional antennas.~However, despite these advantages, unlike the sub-6 GHz band, transmissions in the mmWave band
are highly sensitive to attenuation by the atmosphere, and blockages of communication channels caused by objects such as humans, concrete buildings, and foliage~\cite{blockage1,blockage2}.~In the event of a blockage between a transmitter and the corresponding receiver antenna, the instantaneous transmission rate drops to almost zero~\cite{hashemi2017out}, and hence such blockages increase the end-to-end packet delays in the network~\cite{wang2018joint}.

In this work we consider a single-user network in which a mmWave band downlink communication channel connects the base station to the user. We design efficient scheduling policies which minimize the impact of such blockages by dynamically optimizing a weighted sum of the end-to-end packet delays, and packet transmission costs.~Since blockages are characterized by ``memory'' or temporal correlation, we model the channel using the Gilbert-Elliott model~\cite{sadeghi2008finite}.~Let $s(t)\in \{0,1\}$ denote the state of the channel at time $t$, where $s(t)=0$ denotes that channel is bad, and $s(t)=1$ that the channel is good. The process $\{s(t)\}$ is assumed to be a two-state Markov chain.~Gilbert-Elliott model has been used earlier for modeling mmWave communications~\cite{blockage_model2,blockage_model1,pan2022age}, and allows us to design intelligent scheduling policies which can exploit these temporal correlations of the channel-state process $\{s(t)\}$ in order to generate better estimates of the current state of the channel~\cite{yao2020age}.~Indeed, as is shown in this work, an efficient scheduling policy carries out transmissions only when the conditional probability that the channel is good, conditioned on the information available with the scheduler currently, is sufficiently high, thereby saving the transmission power. This threshold depends upon the queue length. Henceforth, we call this conditional probability as the belief state.
~In contrast with the works such as~\cite{wang2018joint,chen2017delay,zhao2019delay}, which study a seemingly similar setup, we assume that the channel state $s(t)$ is not observed by the transmitter unless a packet transmission is attempted at time $t$, and hence the channel state is ``unobserved''\footnote{Channel state is bad, and blockage are interchangeably used throughout the
paper.}. This is a more realistic assumption since it is not practical for the transmitter to continuously probe the communication channel. Indeed, since the
channel state associated with the mmWave band communication changes at a much rapid rate as compared with that of the sub-6 GHz band, one would need to employ a high frequency probing mechanism in order to track $\{s(t)\}$. This would lead
to a prohibitively high overhead \cite{pomdp_justification1}, which would be exacerbated by the fact that mmWave communication uses large antenna-arrays~\cite{overhead1}.~Though the unobserved channel model is much more realistic, it comes at a significant price since its analysis is quite complex as compared with the observed channel setup due to the following reasons: 
\begin{itemize}
    \item The transmitter now needs to maintain a probabilistic estimate of the channel state, called the belief state. Because this estimate is continuous-valued, it increases the scheduler's memory requirements. Moreover, since scheduling decisions are to be made in real-time, this means that the estimate should also be updated in real-time.
    \item As is shown in our work, the optimal scheduling problem can be posed as a POMDP~\cite{krishnamurthyPOMDP} in which the state of the system is composed of the current queue length, and the belief state. In general, POMDPs are hard to solve~\cite{papadimitriou1987complexity} since the state-space is continuous, and because Dynamic Programming~\cite{bertsekas_book} equations must be solved on an uncountably infinite state-space.
    \item In a real-world setup, underlying system parameters such as the packet arrival rate, and the Markovian transition probabilities associated with the Gilbert-Elliott channel are unknown, so that there is a need for a learning algorithm~\cite{sutton2018reinforcement} that learns an optimal  policy.~However, it is well-known that designing  an efficient learning algorithm for POMDPs is hard~\cite{pmlrcsaba}.
\end{itemize}
Our work overcomes each of the above-mentioned challenges. 
\subsection{Our Contributions}
Our contributions are summarized as follows:
\begin{itemize}
\item For the mmWave network discussed above, we design dynamic scheduling policies that minimize a weighted sum of the average value of the end-to-end packet delays, and packet transmission costs.~In our setup, the transmitter does not need to continually sense the channel, and this reduces the overheads associated with channel sensing. Instead, it gets to observe the channel state only when it attempts a transmission, via ACK\slash NACK mechanism.~Section \ref{sec:Problem-Formulation} shows that the resulting dynamic optimization problem can be formulated as a POMDP in which the system state is composed of (i) the current queue length, and (ii) the belief state.~At each time step, the scheduler uses the current state to decide the number of packets to transmit. 
\item POMDPs are hard to solve in general~\cite{papadimitriou1987complexity}, so that a low-complexity
solution is not guaranteed. However, in Section \ref{sec:Threshold-Policy-Proof} Theorem \ref{th:avg_thresh} we show that the POMDP of interest admits an optimal policy that has
a threshold structure. More specifically, for each value of the queue length $q$, there are multiple non-overlapping intervals $(\tau^{(j)}(q),\tau^{(j+1)}(q))$, $\tau^{(j)}(q)\le\tau^{(j+1)}(q)$, such that it is optimal to transmit $j$ packets when the current queue length is equal to $q$, and the belief state lies within the interval $(\tau^{(j)}(q),\tau^{(j+1)}(q))$. Though previous works on POMDPs derive structural results for POMDPs, they are either restricted to the case of binary action sets \cite{yao2020age,pomdp_twosate_important2} or to the case where the entire system state is described by the belief state~\cite{krishnamurthyPOMDP, krishnamurthy2007structured}.
~To the best of our knowledge, ours is the first work on POMDPs that allows for (i) choosing from a $k$-ary action set with $k>2$, henceforth denoted as a multi-action set, and  
(ii) also has a non-trivial system state, the queue length $Q(t)$, in addition to the belief state associated with the partially observed process.~The key difficulty encountered while going from the binary action setup to the multi-action set case is to show that the decision regions, i.e., the set of belief states in which a given action is optimal for a fixed value of queue length, are contiguous intervals. We show this property in Propositions~\ref{prop:1} and~\ref{prop:3}.~We also show that these decision regions are ``monotonic,'' i.e. keeping the queue length fixed, it is optimal to transmit a higher number of packets only when the belief state is higher.

\item In Section \ref{sec:Simulation-Results} we consider the case when the network parameters such as the Markovian channel parameters and packet arrival probabilities are not known, and the goal is to design an efficient learning algorithm that converges to an efficient scheduling policy. We leverage the threshold structure of an optimal policy that is shown in Theorem \ref{th:avg_thresh}, in order to propose an actor-critic (AC) based learning algorithm. Specifically, we introduce a novel paramerization for generating threshold policies for the scheduling problem, and tune its parameters adaptively using the AC algorithm. We show, by performing extensive simulations, that the performance of the AC algorithm is comparable to that of the optimal policy obtained by using the relative value iteration (RVI) algorithm that has knowledge of the network parameters.

\end{itemize}

\subsection{Related Literature}
We begin by describing works that use a Markov Decision Process (MDP) framework in order to design scheduling policies for mmWave networks.~\cite{yao2019latency} considers a hybrid network architecture in which the scheduler has
to dynamically choose between a mmWave channel that suffers from blockages, and a sub-6 GHz channel so as to minimize the average delay of the system. It derives structural results of an optimal policy under the assumption that the state of the mmWave channel is known by the transmitter at each time. Similar problems are studied in \cite{hashemi2017out,cao2020delay_guidan_similar}. The work~\cite{hashemi2017out} considers the problem of maximizing the network throughput subject to ``delay-like''
constraints when the state of the mmWave channel is known at the transmitter, and derives structural results, while~\cite{cao2020delay_guidan_similar} directly learns the optimal scheduling policy using the Q-learning algorithm \cite{watkins1992q}.~However, these results cannot be used to solve our problem since in our setup the channel state is not observed by the transmitter.
~The works~\cite{chen2017delay,zhao2019delay} consider a single-user network consisting of a fading wireless channel in which the transmitter has to dynamically adjust its transmission power so as to attain an optimal trade-off between packet delays and energy consumption. However, it assumes that the transmitter knows the channel state even without carrying out transmissions, and the channel states are i.i.d. across time.~The work~\cite{wang2018joint} studies a similar problem.~Cross layer design has the potential to minimize traffic latency and also maintain energy efficiency, and hence several works explore this~\cite{collins1999transmission,she2017radio,djonin2007mimo,el2002energy}.~Transmissions in mmWave band use highly directional beams in order to overcome the pathloss suffered by signals. Since these directional beams need to be frequently aligned \cite{gonzalez2017millimeter}, the time associated with this alignment introduces overheads. Thus,~\cite{scheduling_eylium} considers the
problem of minimizing the beam re-alignment time in a multi-user setup,
subject to average throughput constraints and designs an optimal policy that has a psuedo-polynomial time complexity.~The usage of highly directional beams also leads to the problem of intermittent
connectivity in such  networks, and this problem gets exacerbated by the mobility of the receiver.~Typically, in mmWave networks, base stations store packets in cache so that transmissions could be carried out at high rates as soon as communication link is established~\cite{scheduling_cache}.~\cite{scheduling_cache} considers a multi-user mmWave network and designs policies that optimally decide the fraction of cache that is allocated
to data of each user so that the cumulative quality of experience of all the users is maximized.~\cite{singh2015decentralized,singh2018throughput,singh2021adaptive} design policies that maximize the ``timely throughput'', i.e., throughput attained from those data packets that reach the destination node within their deadlines, in unreliable networks.

We now discuss several works which use POMDP framework in order to design efficient algorithms for mmWave networks.~\cite{scheduling1_throughput} designs a policy that maximizes the network throughput, while simultaneously satisfying constraints on the average power consumption.
~\cite{pomdp_d2d} studies the problem of selecting relays in mmWave device-to-device networks, so as to minimize the total number of packets that are lost. The setup therein allows a device to choose between continuing transmissions using the  current relay link, or switch to another relay. The state of the channel connecting two devices is not known globally.~It is shown that the optimal policy has the following structure: it is optimal to switch to a better relay link if the number
of successive ACK failures on the current channel is greater than a certain threshold. Similar
results which address the problem of blockages~\cite{pomdp_APSelction_threshold} have been derived for multiple-access point networks that operate in indoor environments and use mmWave.~The works~\cite{yao2020age,pomdp_twosate_important2,pomdp_whittle1,pomdp_twosate_important,pomdp_twosate_important3,pomdp_aoi_2012} utilize a POMDP framework in order to develop optimal scheduling policies for networks.~Amongst these, \cite{pomdp_whittle1,pomdp_twosate_important,pomdp_twosate_important3} maximize the network throughput, while~\cite{yao2020age, pomdp_twosate_important2,pomdp_aoi_2012}
 minimize the average age-of-information\footnote{Even though \cite{pomdp_aoi_2012} claims that it minimizes delay,
a closer look at equations (7), (8), and (9) of \cite{pomdp_aoi_2012}
reveals that they are minimizing the age of information.} (AoI). To the best of our knowledge, ours is the first work that
minimizes latency in networks for which the channel state is not observed unless a transmission is attempted.

\emph{Notation}: We let $\bN, \bZ_{+}, \bR, \bR_+$ to denote the set of natural numbers, non-negative integers, real numbers, and non-negative real numbers, respectively. We assume that we are equipped with a probability space $(\Omega,\cF,\bP)$.

\section{Problem Formulation\label{sec:Problem-Formulation}}
We present the system model in Section \ref{subsec:System-Model},
and then formulate the optimal scheduling problem as a POMDP in Section \ref{subsec:POMDP_Formulation}.

\subsection{System Model\label{subsec:System-Model}}

\noindent 
\begin{figure}[!t]
\centering
\includegraphics[scale=0.95]{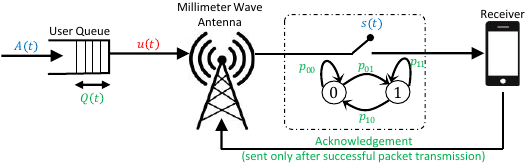} 

\caption{Our system model with the arrival process, $A\left(t\right)$, queue length $Q\left(t\right)$, channel state, $s\left(t\right)$, and the decision variable, $u\left(t\right)$.}
\label{fig:system_model}\vspace{-1em}
 
\end{figure}
Consider a network in which a single user downloads data packets using a downlink communication channel in which the transmitter uses mmWave band for carrying out transmissions (Fig.~\ref{fig:system_model}). Time is discretized into slots $t=0,1,2,\ldots$.~Transmitter stores data packets in a first in first out queue with an infinite buffer size, and fetches packets from it to transmit.~We let $Q\left(t\right)\in\bZ_{+} \label{eq:4.1.6} $ be the
number of packets in the queue at time $t$, and $A\left(t\right)$
be the number of packets that arrive at time $t$. We assume that packet
arrivals $A(t),t\in\bN$ are i.i.d., and uniformly bounded, i.e., $A\left(t\right)\in \mathcal{A}:=\left\{ 0,1,\ldots,M_{a}\right\} $,
where $M_{a}\in\mathbb{N}$. Thus, the distribution of packet arrivals is completely described by the probabilities $p_{i}:=\mathbb{P}(A\left(1\right)=i),~i=0,1,\ldots,M_{a}$.~At each time $t=1,2,\ldots$, the transmitter has to decide the number of packets that are to be attempted for transmission.~We let $u(t)\in \cU := \{0,1,\ldots,M_d\}$ be the number of packets that are attempted for transmission at time $t$. Note that $u(t)=0$ denotes that the transmitter decides to not carry out any transmission at time $t$.~Packet transmissions incur cost, and we let $c(u) \ge 0,~u\in\cU$ be the cost incurred when $u$ packets are attempted transmission. 
\begin{assumption}
    We let $c(0) =0$, and assume that for $u_1,u_2 \in\cU$ satisfying $u_{1}>u_{2}$ we have $c(u_{1}) > c(u_{2})$, i.e., transmitting more packets incurs a higher cost.
\end{assumption}
The line-of-sight between the mmWave transmitter and the receiver is prone to blockages caused due to the presence of objects such as vehicles, humans, buildings, trees~\cite{blockage1,blockage2}. Any transmission attempted when the channel is blocked, fails. The transmitter does not deploy a continual probing mechanism since these incur a lot of overhead \cite{overhead1,overhead2}. Hence, it does not get to observe a blockage unless it transmits packets.~More specifically, the receiver sends an acknowledgement to the transmitter upon receiving packets, so that when a transmission fails, the transmitter does not receive any acknowledgement, and hence concludes that the channel was blocked.~This is the ACK\slash NACK channel feedback \cite{yao2020age,ouyang2015low}, and leads to a ``partially observed'' channel model.~The channel-state process $\{s(t)\}$ is a two-state Markov process~\cite{blockage_model2,blockage_model1}.~$s(t)=0$ denotes that the channel is blocked at time $t$, while $s(t)=1$ denotes that the channel is available. 
~The transition probabilities associated with $s(t)$ are as follows, 
\begin{equation}
\begin{array}{cc}
 & \begin{array}{cc}
s\left(t+1\right)=0\; & s\left(t+1\right)=1\end{array}\\
\begin{array}{c}
s\left(t\right)=0\\
s\left(t\right)=1
\end{array} & \left[\begin{array}{cc}
p_{00}\:\:\:\:\:\:\:\:\:\:\:\:\:\:\:\:\: & p_{01}\:\:\:\\
p_{10}\:\:\:\:\:\:\:\:\:\:\:\:\:\:\:\:\: & p_{11}\:\:\:
\end{array}\right],
\end{array}\label{eq:2.1.1}
\end{equation}
where $p_{00},p_{11}\in(0,1)$. We let $\mu_i, i \in \{0,1\}$ be the stationary probability associated with the channel state process $\{s(t)\}$ being in state $i$, and denote $\mu = (\mu_0,\mu_1)$.~The queue length process $Q(t), t\in\bZ_{+}$ evolves as follows,
\begin{equation}
Q(t+1) = \max\left(0,Q(t)-u(t)s(t)\right)+A(t).\label{eq:2.1.2}
\end{equation}
Define, $\cQ^+(Q,a,d) := \max\left(0,Q-d\right)+a$ so that the queue evolution can be written compactly as $Q(t+1)= \cQ^+(Q(t),A(t),u(t)s(t))$.

\subsection{POMDP Formulation\label{subsec:POMDP_Formulation}}
At each time $t$, the transmitter has to choose the number of packets $u(t)$ that it attempts for transmission, so as to minimize the following objective,
\begin{align}
\underset{T\rightarrow\infty}{\limsup}\:\frac{1}{T}\mathbb{E}\left[\stackrel[t=0]{T-1}{\sum} Q\left(t\right)+ \kappa c(u\left(t\right))\right],\label{def:obj}
\end{align}
where expectation is taken with respect to the probability measure induced by the processes $\{u(t)\},\{A(t)\},\{s(t)\}$, and $\kappa>0$ is a parameter that decides the weightage given to the transmission costs.~The cost in~\eqref{def:obj} is a weighted sum of the following two quantities: (i) average queue length at the transmitter; by Little's Law, this is proportional to the average queueing delay, and (ii) average transmission costs incurred while transmitting data packets. The instantaneous cost in~\eqref{def:obj} captures the tradeoff between energy usage and delay, i.e., transmitting fewer data packets reduces energy usage but increases queueing delay, while transmitting more data packets shortens the queue at the expense of higher energy cost.~Now let $\cF_t$ be the sigma-algebra generated by $\{u(k)\}_{k=1}^{t-1},\{A(k)\}_{k=1}^{t},\{u(k)s(k)\}_{k=1}^{t-1}$. $\cF_t$ denotes the information available with the scheduler at time $t$. A scheduling policy $\pi: \cF_t \mapsto \cU,~t=1,2,\ldots$, decides how many packets should be transmitted at each time $t$, based on the information available with it. Since $\cF_t$ grows with time, it is not convenient to store the entire operational history. Instead,~\eqref{def:obj} can be posed as a POMDP~\cite{krishnamurthyPOMDP,kumar2015stochastic}. The system state for this POMDP at time $t$ is composed of (i) the queue length $Q(t)$, and (ii) the conditional probability of the event that the channel state $s(t)$ is equal to $1$, given the information available with the scheduler. This conditional probability is denoted by $b(t)$, and updated recursively as follows \cite{pomdp_whittle1},
\begin{equation}
b\left(t+1\right)=\begin{cases}
p_{01}& \mbox{ if } u\left(t\right)>0,s\left(t\right)=0,\\
p_{11}& \mbox{ if } u\left(t\right)>0,s\left(t\right)=1,\\
\mathcal{T}\left(b\left(t\right)\right) & \mbox{ if } u\left(t\right)=0,
\end{cases}\label{eq:2.2.2}
\end{equation}
where for $b\in [0,1]$, the function $\mathcal{T}\left(b\right)$ is defined as follows,
\begin{equation}
\mathcal{T}\left(b\right):=b p_{11}+\left(1-b\right)p_{01}.\label{eq:2.2.2.1}
\end{equation}
$b(t)$ is commonly
called the ``belief state'' in the POMDP literature~\cite{krishnamurthyPOMDP,kumar2015stochastic}. We use $\cB$ to denote the space in which $b(t)$ lives; $\cB$ is discussed below in detail. It can be shown that $b(t)$ serves as a \textit{sufficient statistic}~\cite[Appendix A]{sufficient_statistics} for the purpose of minimizing the cost~\eqref{def:obj}~\cite{krishnamurthyPOMDP}. Thus, we can replace the original system by one in which the system state is taken to be $(Q(t),b(t))$, the instantaneous cost incurred at time $t$ is equal to $Q(t)+ \kappa c(u(t))$, the evolution of the state $(Q(t),b(t))$ is described by~\eqref{eq:2.1.2} and~\eqref{eq:2.2.2}, and any transmission attempt at time $t$ is successful with a probability $b(t)$. We denote the combined state-space in which $(Q(t),b(t))$ resides by $\cX := \bZ_{+} \times \cB$.~With this setup in place, a scheduling policy $\pi$
for each time $t=1,2,\ldots$, chooses $u(t)$ as a function of the state $(Q(t),b(t))$. The objective is to solve the following POMDP:
\begin{align}
\min_{\pi} \underset{T\rightarrow\infty}{\limsup}\:\frac{1}{T}\mathbb{E}_{\pi}\left[\stackrel[t=0]{T-1}{\sum} Q\left(t\right)+ \kappa c(u\left(t\right))\right].\label{def:pomdp}
\end{align}

\textbf{Belief Space $\mathcal{B}$}: We now briefly describe the space in which the belief state $b(t)$ lives.~If the initial belief state is $b$, and then no transmission occurs for $k$ consecutive time steps, then the resulting belief state is given by $\mathcal{T}^{k}\left(b\right)$, 
\[
\mathcal{T}^{k}\left(b\right):=\underset{k-times}{\underbrace{\mathcal{T}\circ\mathcal{T}\cdots\circ\mathcal{T}}}\left(b\right)\;,\:k=1,2,\ldots,
\]
where $\mathcal{T}^{0}\left(b\right):=b$, for two functions $f,g$, $f\circ g$ denotes their composition, and $\mathcal{T}\left(\cdot\right)$ is given
by (\ref{eq:2.2.2.1}).~Define the following sets, 
\begin{eqnarray}
\mathcal{B}_{0} & := & \left\{ \mathcal{T}^{k}\left(b\left(0\right)\right)\:|\:k=0,1,2,\ldots\right\} ,\label{eq:2.2.0.1}\\
\mathcal{B}_{1} & := & \left\{ \mathcal{T}^{k}\left(p_{01}\right)\:|\:k=0,1,2,\ldots\right\} ,\label{eq:2.2.0.2}\\
\mathcal{B}_{2} & := & \left\{ \mathcal{T}^{k}\left(p_{11}\right)\:|\:k=0,1,2,\ldots\right\},\label{eq:2.2.0.3}\\
\mbox{ and }\mathcal{B} & := & \mathcal{B}_{0}\cup\mathcal{B}_{1}\cup\mathcal{B}_{2}\label{eq:2.2.0.0}.
\end{eqnarray}
Let $t_{0}:=\min\left\{ t\geq0:u\left(t\right)>0\right\} $, denote the first time step when a transmission occurs. Then, for $t\leq t_{0}$ we have $b(t)\in\mathcal{B}_{0}$. Since $u\left(t_{0}\right)>0$, it follows from (\ref{eq:2.2.2}) that $b\left(t_{0}+1\right)\in\left\{ p_{01},p_{11}\right\} $. Hence, for $t>t_{0}$ we have $b(t)\in\mathcal{B}_{1}\cup\mathcal{B}_{2}$.
\section{Threshold Structure of Optimal Policy}
\label{sec:Threshold-Policy-Proof}
In general, solving a POMDP is hard. For example,~\cite{papadimitriou1987complexity} shows that in general a POMDP is PSPACE hard.~However, it turns out that many POMDPs of practical interest admit optimal policies with an appealing structure. If one is able to show such a structural result, then it can be used to obtain
a low-complexity implementation of the optimal policy.~This also makes it easier to search for an optimal policy since one can now restrict the search space.~We begin by defining what is meant by threshold structure for a policy for the POMDP~\eqref{def:pomdp}, and then show in this section that~\eqref{def:pomdp} admits an optimal policy that has this structure. 
\begin{definition}
We say that a policy $\pi: \mathcal{X} \mapsto \{0,1,\ldots,M_d\}$ is of threshold-type if for each $q\in \bZ_{+}$, there are threshold values $\tau^{(j)}(q),~j=0,1,\ldots,M_d +1$ with $\tau^{(0)}(q)=0$, $\tau^{(M_{d}+1)}(q)=1$, and $\tau^{(j)}(q) \leq \tau^{(j+1)}(q)$ for $j=0,1,\ldots,M_d$, such that the decision chosen when $Q(t)=q$ can be described as follows:
\begin{itemize}
    \item if $b(t)\in \left[\tau^{(j)}(q),\tau^{(j+1)}(q)\right)$, where $j\in\{0,1,\ldots,M_d -1\}$, then it transmits $j$ packets.
    \item if $b(t)\in \left[\tau^{(M_d)}(q),\tau^{(M_{d+1})}(q)\right]$, then it transmits $M_d$ packets.
\end{itemize}
\end{definition}
We note that in order to describe a general stationary deterministic policy, we need to mention the decision $u$ for each possible $(q,b)$. On the other hand, to describe a threshold policy, we simply need to provide $M_d + 1$ thresholds $\{\tau^{(j)}(q)\}_{j=0}^{M_d}$, and thus this structure yields a vast simplification. 

Next, we begin by studying the $\beta$-discounted cost problem in Section \ref{subsec:General-threshold} since it is simpler than the original average cost POMDP (\ref{def:pomdp}). We derive structural results for the $\beta$-discounted cost problem.~Section~\ref{subsec:Optimal-Policy-Average} extends this result to the original average-cost problem~(\ref{def:pomdp}).~Section~\ref{subsec:Optimal-Policy-Average-Computing} discusses how to compute an optimal policy.

\subsection{Discounted POMDP}
\label{subsec:General-threshold}

\noindent 
\begin{figure}[t]
\begin{centering}
\includegraphics[scale=0.9]{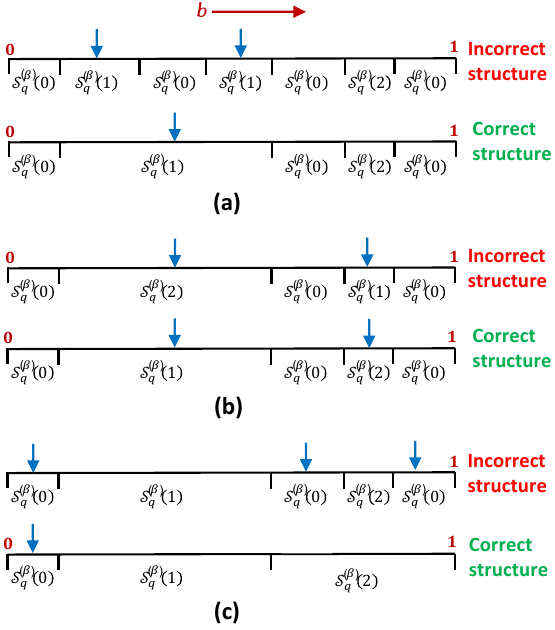} 
\par\end{centering}
\caption{Pictorial depiction of the threshold structure of policy that is optimal for the $\beta$-discounted problem for $M_{d}=2$: (a) Diagram showing that the set $\mathcal{S}_{q}\ubeta\left(u\right)$
are contiguous intervals for $u\protect\geq1$. The blue arrows in
``incorrect structure'' show that the set $\mathcal{S}_{q}\ubeta\left(1\right)$
consists of two disjoint intervals. This is not possible. (b) Diagram
showing the ordering of the set $\mathcal{S}_{q}\ubeta\left(u\right)$
for $u\protect\geq1$. The blue arrows in ``incorrect structure''
show that the elements of set $\mathcal{S}_{q}\ubeta\left(2\right)$
is smaller than that of $\mathcal{S}_{q}\ubeta\left(1\right)$.
This is not possible. (c) (a) Diagram showing that the set $\mathcal{S}_{q}\ubeta\left(0\right)$
is a contiguous interval. The blue arrows in ``incorrect structure''
show that the set $\mathcal{S}_{q}\ubeta\left(0\right)$ consists
of two disjoint intervals. This is not possible.}
\label{fig:threshold_diag}\vspace{-1em}
 
\end{figure}
We will analyze the following $\beta$-discounted POMDP
\begin{align}
\min_{\pi}\mathbb{E}_{\pi}\left[\stackrel[t=0]{\infty}{\sum}\beta^{t}\left(Q\left(t\right)+ \kappa c(u\left(t\right))\right)\right],\label{def:discountedpomdp}
\end{align}
where $\beta \in (0,1)$ is the discount factor.~Since the state-space $\cX$ associated with (\ref{def:discountedpomdp}) is infinite, it may not admit an optimal policy 
~\cite[p. 52]{hernandez2012discrete}.~However, it is shown in Appendix~\ref{app:stat_optimal} that there exists an optimal policy that solves (\ref{def:discountedpomdp}).~Moreover, it is also shown that for the purpose of solving (\ref{def:discountedpomdp}), it suffices to restrict to the class of stationary policies.~We now describe the dynamic programming equation \cite{bertsekas_book} corresponding to the problem~\eqref{def:discountedpomdp}. Let $J\ubeta: \cX \mapsto \mathbb{R}$ be the value function of (\ref{def:discountedpomdp}), i.e., $J\ubeta\left(q,b\right)$ is the optimal discounted cost~\eqref{def:discountedpomdp} incurred by the system when it starts in the state $\left(q,b\right)$.~Let $\tilde{J}\ubeta(q,b;u)$ be the cumulative discounted
cost incurred by the system when it starts in state $(q,b)$, and action $u$ is applied at time $t=0$, following which an optimal policy is implemented from time $t=1$ onwards, i.e., 
\al{
& \tilde{J}\ubeta\left(q,b;0\right) : =q+\beta\stackrel[i=0]{M_{a}}{\sum}p_{i}J\ubeta\left(\cQ^+(q,i,0),\mathcal{T}\left(b\right)\right),\notag \\
& \tilde{J}\ubeta\left(q,b;u\right) : =q+ \kappa c(u) \notag\\
& +\beta\stackrel[i=0]{M_{a}}{\sum}p_{i}\left(b J\ubeta\left(\cQ^+(q,i,u),p_{11}\right) + (1-b)\right. \notag \\
& \times \left. J\ubeta\left(\cQ^+(q,i,0),p_{01}\right)\right),~u\in\{1,2,\ldots,M_d\}. \label{eq:3.3}
}
It is shown in Appendix \ref{app:Boundedness-Proof} Lemma~\ref{lem:SennottC1} that $J\ubeta\left(q,b\right)$ is finite for all $\beta\in(0,1)$ and $(q,b)\in \cX$. Hence, from~\cite[Proposition 1]{sennott1989average} we conclude that the value function $J\ubeta\left(\cdot,\cdot\right)$ of the discounted cost problem (\ref{def:discountedpomdp}) solves the dynamic programming equation, i.e., we have,
\begin{equation}
J\ubeta\left(q,b\right)=\underset{u\in\cU }{\min}~\tilde{J}\ubeta\left(q,b;u\right).\label{eq:3.1}
\end{equation}

The value function $J\ubeta(\cdot,\cdot)$ can be obtained by using the value iteration algorithm \cite[Proposition 3]{sennott1989average}, which we describe now. Let $J\ubeta_{n}\left(q,b\right)$ be the estimate of the value function generated by the value iteration algorithm during the $n$-th iteration. During the $n$-th iteration, it first computes the quantities $\tilde{J}\ubeta_{n}\left(q,b;u\right)$ as follows,
\begin{align}
&\tilde{J}\ubeta_{n}\left(q,b;0\right) =q+\beta\stackrel[i=0]{M_{a}}{\sum}p_{i}J\ubeta_{n-1}\left(\cQ^+(q,i,0), \mathcal{T}\left(b\right)\right),\label{eq:3.2.1} \\
& \tilde{J}\ubeta_{n}\left(q,b;u\right)=q+ \kappa c(u) \notag \\
& +\beta\stackrel[i=0]{M_{a}}{\sum}p_{i}\left(b J\ubeta_{n-1}\left(\cQ^+(q,i,u),p_{11}\right) + (1-b)\right. \notag\\
&\left.\times J\ubeta_{n-1}\left(\cQ^+(q,i,0),p_{01}\right)\right),~u\in\{1,2,\ldots,M_{d}\}.\label{eq:3.3.1}
\end{align}
Finally, $J\ubeta_n(q,b)$ gets updated as follows,
\begin{equation}
J_{n}\ubeta\left(q,b\right)=\underset{u\in\cU }{\min}~\tilde{J}\ubeta_{n}\left(q,b;u\right).\label{eq:3.0.2}
\end{equation}
We initialize $J_{0}\ubeta\left(q,b\right)=0,~\forall\left(q,b\right)\in\mathcal{X}$. The iterates $J_{n}\ubeta\left(q,b\right)$ converge to the value function~\cite[Proposition 3]{sennott1989average}, i.e., 
\al{
\underset{n\rightarrow\infty}{\lim}J_{n}\ubeta\left(q,b\right) = J\ubeta\left(q,b\right), \forall\left(q,b\right)\in\cX. \label{eq:Jn_convg_J}
}
We now present the main result of this section.
\begin{theorem}
\label{thm:DiscountedThreshold}
Consider the $\beta$-discounted POMDP (\ref{def:discountedpomdp}). There exists a threshold policy that is optimal. 
\end{theorem}
\begin{IEEEproof}
Let $\mathcal{S}_{q}\ubeta\left(u\right)$ denote the set of belief
states $b$, in which it is optimal to transmit $u$ packets when the
current queue length is equal to $q$, i.e., 
\[
\mathcal{S}_{q}\ubeta\left(u\right):=\left\{ b:\tilde{J}\ubeta\left(q,b;u\right)= J\ubeta\left(q,b\right)\right\}. 
\]
The proof follows from the following three propositions.

\begin{proposition}\label{prop:1}
For each $q\in \bZ_{+}$, $\beta\in (0,1)$, and $u \in \{1,2,\ldots,M_d\}$, the sets $\mathcal{S}_{q}\ubeta\left(u\right)$ are contiguous intervals~(Fig. \ref{fig:threshold_diag}a).
\end{proposition}
\proof  Consider $b_{1},b_{2}\in\mathcal{S}_{q}\ubeta\left(u\right)$.
We will show that $\gamma b_{1}+\left(1-\gamma\right)b_{2}\in\mathcal{S}_{q}\ubeta\left(u\right)$ for all $\gamma\in\left[0,1\right]$. Let $\tilde{b}$ denote $\gamma b_{1}+\left(1-\gamma\right)b_{2}$. We will show $\tilde{J}\ubeta\left(q,\tilde{b};u\right) = J\ubeta\left(q,\tilde{b}\right)$, which will prove the claim.~We have, 
\begin{align*}
J\ubeta\left(q,\tilde{b}\right) & \geq \gamma J\ubeta\left(q,b_{1}\right)+\left(1-\gamma\right)J\ubeta\left(q,b_{2}\right)\nonumber \\
 & =  \gamma \tilde{J}\ubeta\left(q,b_{1};u\right)+\left(1-\gamma\right)\tilde{J}\ubeta\left(q,b_{2};u\right)\\
 & =  q+\kappa c(u)+\beta\stackrel[i=0]{M_{a}}{\sum}p_{i}\left(\tilde{b}J\ubeta\left(\cQ^+(q,i,u),p_{11}\right)\right.\\
 &   \quad\quad\quad\quad\left.\vphantom{\beta\left(\tilde{b}J\ubeta\left(\left(q+i-j\right)_{0}^{M_{q}},p_{11}\right)+\right.}+\left(1-\tilde{b}\right)J\ubeta\left(\cQ^+(q,i,0),p_{01}\right)\right)\\
 & =  \tilde{J}\ubeta\left(q,\tilde{b};u\right)\\
 & \geq  J\ubeta\left(q,\tilde{b}\right),
\end{align*}
where the first inequality follows from Lemma \ref{lem:lemma3}, the second and third equalities follow from (\ref{eq:3.3}), and the last inequality follows from (\ref{eq:3.1}).
This gives us $J\ubeta\left(q,\tilde{b}\right)\geq \tilde{J}\ubeta\left(q,\tilde{b};u\right)\geq J\ubeta\left(q,\tilde{b}\right)$.
The proof is then completed by observing that for the above set of inequalities to hold,
we must have $\tilde{J}\ubeta\left(q,\tilde{b};u\right) = J\ubeta\left(q,\tilde{b}\right)$.

\begin{proposition}\label{prop:2}
Fix a discount factor $\beta\in(0,1)$, and consider $u_{1},u_{2}\in\{1,2,\ldots,M_{d}\}$ satisfying $u_{1}<u_{2}$. For each $q\in \bZ_{+}$, if $b_{1}\in\mathcal{S}_{q}\ubeta\left(u_{1}\right)$,~$b_{2}\in\mathcal{S}_{q}\ubeta\left(u_{2}\right)$,
 then we must have $b_{1}\leq b_{2}$. Thus, it is optimal to transmit a larger number of packets only when the channel is less likely to be blocked (Fig.~\ref{fig:threshold_diag}b).   
\end{proposition} 

\proof We prove this by contradiction. On the contrary, let us assume that $b_{2}<b_{1}$. We then have that, 
\begin{align}
 &  \tilde{J}\ubeta\left(q,b;u_{2}\right)-\tilde{J}\ubeta\left(q,b;u_{1}\right)\nonumber \\
 & =  \kappa\left(c(u_{2})-c(u_{1})\right)+\beta b\stackrel[i=0]{M_{a}}{\sum}p_{i}\left(J\ubeta\left(\cQ^+(q,i,u_2),p_{11}\right)\right.\nonumber \\
 &  \quad\quad\quad\quad\quad\left.\vphantom{\beta b\left(J\ubeta\left(\cQ^+(q,i,u_2),p_{11}\right)-\right.}-J\ubeta\left(\cQ^+(q,i,u_1),p_{11}\right)\right).\label{eq:Th6.1}
\end{align}
Now since $u_{1}<u_{2}$, and $J\ubeta\left(q,b\right)$ is monotonic non-decreasing (Lemma \ref{lem:lemma4-1}), we have that $J\ubeta\left(\cQ^+(q,i,u_2),p_{11}\right) - J\ubeta\left(\cQ^+(q,i,u_1),p_{11}\right)\le 0$. Thus, $J\ubeta\left(q,b;u_{2}\right)-J\ubeta\left(q,b;u_{1}\right)$
is monotonic non-increasing function of $b$. Since $u_1$ is optimal for $b_{1}$, we have $J\ubeta\left(q,b_{1};u_{1}\right)\leq J\ubeta\left(q,b_{1};u_{2}\right)$. However, since $J\ubeta\left(q,b;u_{2}\right)-J\ubeta\left(q,b;u_{1}\right)$
is a monotonically non-increasing function of $b$, we have $\tilde{J}\ubeta\left(q,b_{2};u_{1}\right)\leq \tilde{J}\ubeta\left(q,b_{2};u_{2}\right)$. But then, if $\tilde{J}\ubeta\left(q,b_{2};u_{1}\right)\leq \tilde{J}\ubeta\left(q,b_{2};u_{2}\right)$ holds, then we conclude that $u_2$ is optimal for $b_{2}$, i.e. $b_{2}\notin\mathcal{S}_{q}\ubeta\left(u_{2}\right)$.~However, this leads to a contradiction since we assumed that $b_{2}\in\mathcal{S}_{q}\ubeta\left(u_{2}\right)$. Hence, we must have $b_1\le b_2$. This completes the proof.

\begin{proposition}\label{prop:3}
For each $q\in \bZ_{+}$, and $\beta\in(0,1)$, $\mathcal{S}_{q}\ubeta\left(0\right)$ is a contiguous
interval (Fig.~\ref{fig:threshold_diag}c).
\end{proposition} 
\proof Fix a $u\in\{1,2,\ldots,M_d\}$. We will prove that each element of $\mathcal{S}_{q}\ubeta\left(0\right)$
is less than each element of $\mathcal{S}_{q}\ubeta\left(u\right)$. Consider $b_{1}\in\mathcal{S}_{q}\ubeta\left(0\right)$,
$b_{2}\in\mathcal{S}_{q}\ubeta\left(u\right)$,~and assume on the contrary that $b_{2}<b_{1}$. We have the following from Lemma \ref{lem:lemma3},
\al{
\tilde{J}\ubeta\left(q,b;0\right) & = q+\beta\stackrel[i=0]{M_{a}}{\sum}p_{i}J\ubeta\left(\cQ^+(q,i,0),\mathcal{T}\left(b\right)\right)\nonumber \\
 & \geq \tilde{J}\ubeta\left(q,b;\widetilde{0}\right), \label{eq:Th7.1}
}
where,
\nal{
\tilde{J}\ubeta\left(q,b;\widetilde{0}\right) & : = q+\beta\stackrel[i=0]{M_{a}}{\sum}p_{i}\left(b J\ubeta\left(\cQ^+(q,i,0),p_{11}\right)\right.\\
 & \left.+\left(1-b\right)J\ubeta\left(\cQ^+(q,i,0),p_{01}\right)\right).
}

Now, similar to (\ref{eq:Th6.1}), we can show that $\tilde{J}\ubeta\left(q,b;u\right)-\tilde{J}\ubeta\left(q,b;\widetilde{0}\right)$
is a monotonically non-increasing function of $b$. Since $0$ is an optimal action in state $(q,b_1)$, we have $\tilde{J}\ubeta\left(q,b_{1};0\right)\leq \tilde{J}\ubeta\left(q,b_{1};u\right)$.
Hence, from (\ref{eq:Th7.1}) we have that $\tilde{J}\ubeta\left(q,b_{1};\widetilde{0}\right)\leq \tilde{J}\ubeta\left(q,b_{1};u\right).$
But then since $\tilde{J}\ubeta\left(q,b;u\right)-\tilde{J}\ubeta\left(q,b;\widetilde{0}\right)$
is monotonically non-increasing in $b$, this shows us that $\tilde{J}\ubeta\left(q,b_{2};\widetilde{0}\right)\leq \tilde{J}\ubeta\left(q,b_{2};u\right)$.
From this, we conclude that $u$ is not optimal for $b_2$, which is a contradiction since we assumed that $b_{2}\in\mathcal{S}_{q}\ubeta\left(u\right)$. Hence, we must have $b_1 \le b_2$. This completes the proof.
\end{IEEEproof}

\subsection{Optimal Policy for the Average Cost Problem\label{subsec:Optimal-Policy-Average}}
This section extends the structural result of Section \ref{subsec:General-threshold} to the average cost problem (\ref{def:pomdp}). We make the following assumption on the underlying POMDP while analyzing the average cost problem.
\begin{assumption} \label{assumption2} 
The packet arrival probabilities satisfy $p_{i}>0$ for all $i=1,2,\ldots,M_d$, and the Markovian channel parameters satisfy $0<p_{01},p_{11}<1$. 
\end{assumption}
\begin{assumption} \label{assum:stability}
    The system satisfies the stability assumption $M_d \mu_1 > \bE[A(1)]$, where $\mu_1$ is the stationary probability of the Markovian channel being in state $1$, $A(t)$ is the number of packets that arrive at time $t$, and $M_d$ is the maximum number of packets that can be attempted for transmission by the transmitter at time $t$.
\end{assumption}
The following result establishes a connection between the average
cost POMDP (\ref{def:pomdp}) and the $\beta$-discounted cost POMDP (\ref{def:discountedpomdp}) as the discount factor $\beta$ approaches one. 
\begin{proposition}
\label{thm:SennotsTheorem}Consider the discounted cost POMDP (\ref{def:discountedpomdp}) and the average cost POMDP (\ref{def:pomdp}).~Let Assumptions \ref{assumption2} and \ref{assum:stability} hold. Let $\{\beta_n\}_{n\in\bN}$ be a sequence of discount factors such that $\beta_n \in (0,1)$ and $\beta_n \to 1$. Let $\{\pi_{\beta_{n}}^{\star}\}_{n\in\bN}$ be the corresponding sequence of optimal policies, i.e. $\pi_{\beta_{n}}^{\star}$ is optimal for (\ref{def:discountedpomdp}) with the discount factor set equal to $\beta_n$. Then, there exists a subsequence of $\{\pi_{\beta_{n}}^{\star}\}$, say $ \pi_{\beta_{n_k}}^{\star}, k=1,2,\ldots$, and a stationary policy $\pi^{\star}$ such that $\pi^{\star}$ is a limit point of $\left\{ \pi_{\beta_{n_k}}^{\star}\right\}$. Moreover, $\pi^{\star}$ is an optimal policy for the average cost POMDP (\ref{def:pomdp}).
\end{proposition}
\begin{IEEEproof}
The proof follows from the lemma in \cite[p.~3]{sennott1989average}, once we have shown that assumptions $1$, $2$, $3$, and $3^{\star}$ made therein hold true for our setup. Lemmas \ref{lem:SennottC1}, \ref{lem:SennottC2}, \ref{lem:SennottC3}, and \ref{lem:SennottC3star} in Appendix \ref{app:Boundedness-Proof} validate assumptions $1$, $2$, $3$, and $3^{\star}$ respectively. This concludes the proof.
\end{IEEEproof}
\begin{theorem}\label{th:avg_thresh}
Under Assumptions~\ref{assumption2} and~\ref{assum:stability}, for the average cost POMDP~\eqref{def:pomdp}, there exists an optimal policy that has a threshold structure.
\end{theorem}
\begin{IEEEproof}
Consider a sequence $\{\beta_n\}_{n\in\bN}, \beta_n \in (0,1)$ of discount factors converging to $1$.~It follows from Proposition \ref{thm:SennotsTheorem} that there exists a sequence of policies $ \pi_{\beta_{n_k}}^{\star}, k=1,2,\ldots$, where $\pi_{\beta_{n_k}}^{\star}$ is optimal for the $\beta_{n_k}$-discounted POMDP, such that there is a limiting policy $\pi^{\star}$ which is also optimal for the average cost POMDP. It follows from Theorem \ref{thm:DiscountedThreshold} that each such policy $\pi_{\beta_{n_k}}^{\star}$ has a threshold structure. The proof is then completed by observing that any limit of a sequence of threshold-type policies is also of threshold-type.
\end{IEEEproof}
\subsection{Computing Optimal Policy }\label{subsec:Optimal-Policy-Average-Computing}
We now show that an optimal policy for the average cost POMDP (\ref{def:pomdp}) can be obtained by using the Relative Value Iteration (RVI) algorithm~\cite[Ch: 4]{bertsekas_book}. RVI solves the following average
cost optimality equation (ACOE) \cite{arapostathis1993discrete},
\begin{alignat}{1}
\zeta\ust +h\left(q,b\right) & =\underset{u\in\cU }{\min}~\tilde{h}\left(q,b;u\right),~\forall (q,b)\in \cX,\label{eq:4.1.1.0}
\end{alignat}
where,
\begin{equation}
\tilde{h}\left(q,b;0\right):=q+\stackrel[i=0]{M_{a}}{\sum}p_{i}h\left(\cQ^+(q,i,0),\mathcal{T}\left(b\right)\right),\label{eq:4.1.1.a}
\end{equation}
and 
\begin{align}
& \tilde{h}\left(q,b;u\right):=q+ \kappa c(u)+\stackrel[i=0]{M_{a}}{\sum}p_{i}\left(b h\left(\cQ^+(q,i,u),p_{11}\right)\right.\notag\\
& \left.\vphantom{\left(b h\left(\left(q+i-j\right)_{0}^{M_{q}},p_{11}\right)+\right.}+\left(1-b\right)h\left(\cQ^+(q,i,0),p_{01}\right)\right),~u \in \{1,2,\ldots,M_{d}\}. \label{eq:4.1.1.b}
\end{align}
Note that $h$ is the relative value function \cite[Chapter 4]{bertsekas_book} and $\zeta\ust$ is the average cost.~Before proceeding to solve (\ref{eq:4.1.1.0}), we firstly show that (\ref{eq:4.1.1.0}) has a bounded solution \cite{derman1966denumerable}.

\begin{proposition}
\label{thm:RossTheorem}Consider the average cost POMDP (\ref{def:pomdp})
and let Assumption \ref{assumption2} hold. Then, there is a bounded $h: \cX \mapsto \bR$ and a $\zeta\ust \in\mathbb{R}$ that solves (\ref{eq:4.1.1.0}). Also, $\zeta\ust$
is the optimal average cost of (\ref{def:pomdp}) that is independent of the
initial state $(Q(0),b(0))$.
\end{proposition}
\begin{IEEEproof}
The proof follows from Lemmas \ref{lem:SennottC1}-\ref{lem:SennottC3star} and the theorem in \cite[p.~628]{sennott1989average}. 
\end{IEEEproof}
We now discuss how to solve~(\ref{eq:4.1.1.0}) using the RVI algorithm.~Let $\left(q_{o},b_{o}\right)$ be a ``reference state'' that is chosen by the scheduler. RVI proceeds
as follows. It keeps an estimate of the relative values $h(q,b)$, and updates them iteratively. Let $h_n(q,b)$ denote these estimates during the $n$-th iteration.~In the $n$-th iteration, it firstly computes $\tilde{h}_{n}(q,b;0)$ using the quantities $h_{n-1}(q,b)$ as follows,
\begin{align}
& \tilde{h}_{n}\left(q,b;0\right)=q+\stackrel[i=0]{M_{a}}{\sum}p_{i}h_{n-1}\left(\cQ^+(q,i,0),\mathcal{T}\left(b\right)\right),\label{eq:4.1.2.a} \\
&\tilde{h}_{n}\left(q,b;u\right)=q+ \kappa c(u)+\stackrel[i=0]{M_{a}}{\sum}p_{i}\left(b h_{n-1}\left(\cQ^+(q,i,u),p_{11}\right)\right.\notag\\
&\left. +\left(1-b\right)h_{n-1}\left(\cQ^+(q,i,0),p_{01}\right)\right),~u\in \{1,2,\ldots,M_{d}\}.\notag
\end{align}

Next, it computes the quantities $\tilde{h}\ust_{n}(q,b)$ as follows

\begin{align}
\tilde{h}_{n}\ust\left(q,b\right)  :=\underset{u\in\cU }{\min}~\tilde{h}_{n}\left(q,b;u\right).\label{eq:4.1.2}
\end{align}
Finally, $h_n(q,b)$ gets updated as follows,
\begin{align}
    h_{n}\left(q,b\right)  =\tilde{h}\ust_{n}\left(q,b\right)-\tilde{h}\ust_{n}\left(q_{o},b_{o}\right).\label{eq:4.1.3}
\end{align}
Iterations are performed until the quantity $|h_{n}\left(q,b\right)-h_{n-1}\left(q,b\right)|$ becomes less than a certain threshold for all $\left(q,b\right)$.

\section{Simulation Results\label{sec:Simulation-Results}}
To compute an optimal policy we need to solve the RVI equations (\ref{eq:4.1.2.a})-(\ref{eq:4.1.3}). Since the state space $\cX=Z_+ \times \mathcal{B}$ is infinite, we truncate it as follows while implementing RVI.~The sets $\mathcal{B}_{0}$, 
 $\mathcal{B}_{1}$, and $\mathcal{B}_{2}$ are truncated by restricting the value of $k$ in (\ref{eq:2.2.0.1})-(\ref{eq:2.2.0.3}) to the set $\left\{ 0,1,\ldots,10\right\}$.~We also truncate the queue length to $10$ packets.~The cost of transmitting $u$ packets at each time is taken to be $c\left(u\right)= \exp(u) - 1$, and unless specified we set $\kappa =1$ in the objective of POMDP (\ref{def:pomdp}).~Since most reinforcement learning (RL) algorithms are designed for maximizing average reward rather than minimizing average cost, hence within this section we pose our problem as a reward maximization problem, one in which the instantaneous reward function is taken to be negative of the cost function, i.e. $(10 + \kappa c(M_d))- (Q + \kappa c(u))$. This ensures compatibility with standard RL frameworks without altering the optimal policy.

\begin{figure}[t]
	\begin{centering}
		\includegraphics[scale=0.34]{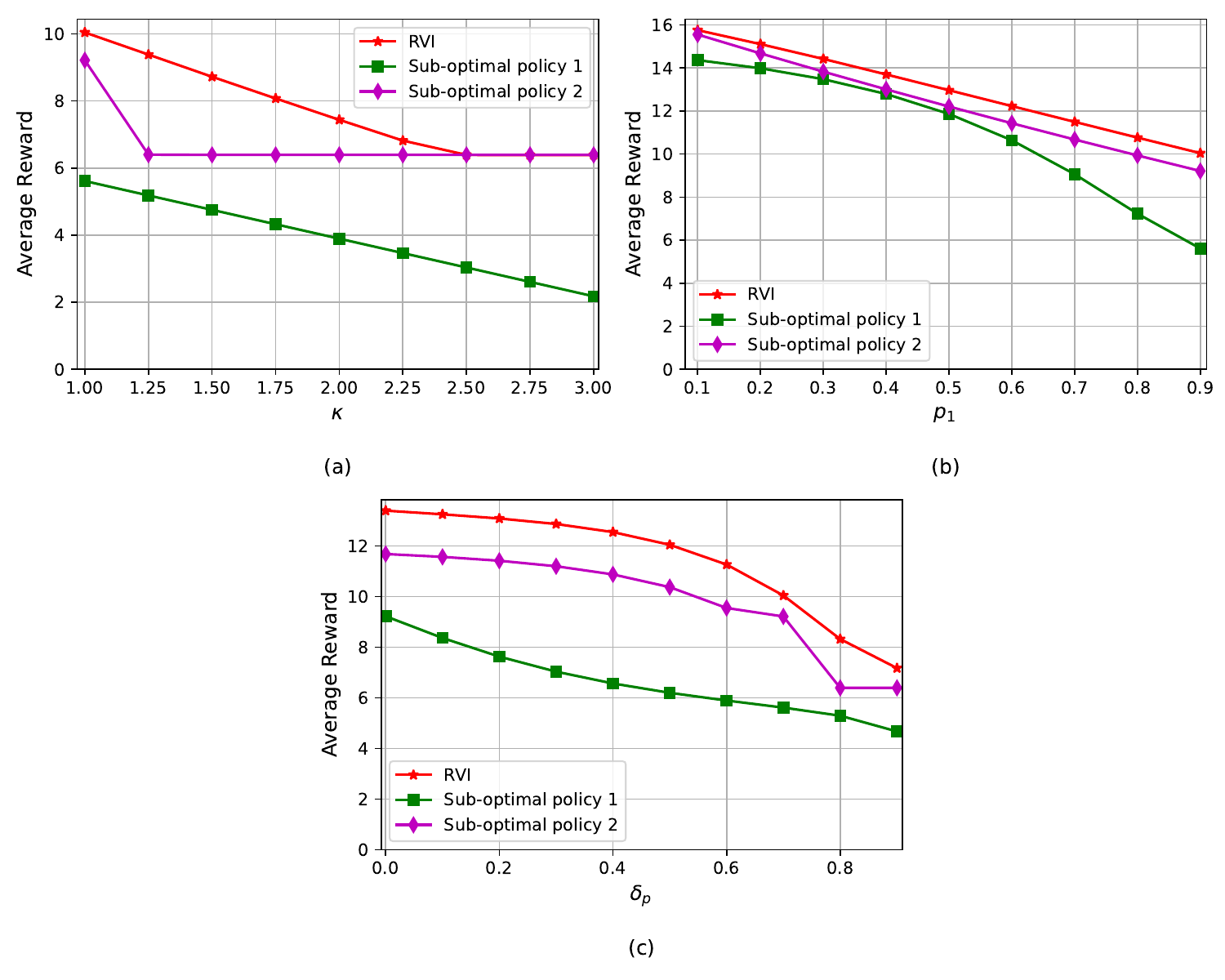} 
		\par\end{centering}
	\caption{Comparison of optimal policy obtained using RVI with two sub-optimal policies: (a) The weight $\kappa$ in~\eqref{def:pomdp} is varied with fixed $p_{01} = 0.2$, $p_{11} = 0.9$, $p_1 = 0.9$, (b) Arrival probability $p_{1}$ is varied with fixed $\kappa =1$, $p_{01} = 0.2$, $p_{11} = 0.9$, and (c) Markovian blockage parameters are varied, here $\delta_p = |p_{11} - p_{01}|$ with $p_{11}$ fixed at $0.9$, as $\delta_p$ is varied and $\kappa =1$.}
	\label{fig:sub_opt_policies}
	
\end{figure}

\subsection{Known network parameters}
Unless mentioned otherwise, the network parameters are chosen as follows. The Markovian channel transition probabilities are set to $p_{01}=.2$, $p_{11}=.9$. $M_d$, which is the maximum number of packets that can be transmitted at each time is set to $2$, and $M_a$ the maximum number of packets that can arrive at each time is taken to be $1$. Arrival probability $p_1$ is taken to be $.9$. We compare the optimal policy for the average cost POMDP (\ref{def:pomdp}) obtained using the RVI as described in Section \ref{subsec:Optimal-Policy-Average-Computing} with the following two policies. 
\begin{itemize}
\item \emph{Sub-optimal policy 1}: Always transmits a single packet at each time step irrespective of the current value of the system state. 
\item \emph{Sub-optimal policy 2}: Assumes that the channel state is i.i.d. across time with the probability of channel in good state equal to the stationary probability that the original Gilbert-Elliott channel is in good state. It then solves an average cost MDP for this i.i.d. channel using the RVI to get this sub-optimal policy.
\end{itemize}
Fig.~\ref{fig:sub_opt_policies} compares the performance of the optimal policy with the above described sub-optimal policies, as $\kappa, p_{1}$, and the quantity $|p_{11} - p_{01}|$ is varied. Note that $|p_{11} - p_{01}|$ is a measure of how ``Markovian'' the channel is. A higher value of $|p_{11} - p_{01}|$ indicates that the channel is ``more Markovian.''~The optimal policy is seen to outperform the suboptimal policies. Specifically, Fig.~\ref{fig:sub_opt_policies}(a) shows that as $\kappa$ is increased, the average rewards decrease. Fig.~\ref{fig:sub_opt_policies}(b) shows that the average reward decreases as $p_{1}$ is increased. This is because due to a higher value of $p_{1}$, the queue length increases, which leads to higher instantaneous costs and more frequent transmissions in order to manage the queue. Both of these result in a decrease in the average reward. Fig.~\ref{fig:sub_opt_policies}(c) depicts the variations with respect to $|p_{11} - p_{01}|$.~Fig.~\ref{fig:sub_opt_policies}(c) shows that as the channel becomes more Markovian, the gap in the performance of the optimal policy as compared with the sub-optimal policies increases. This is because none of the two sub-optimal policies account for the Markovian nature of the channel.
\subsection{Actor-Critic Algorithm} 
In this section we use the actor-critic (AC) algorithm~\cite{konda1999actor} for learning an optimal policy. Since it is a model-free learning algorithm, it doesn't need to maintain an estimate of the network parameters. Since the AC does not solve RVI, within this section we do not truncate the belief-space, however queue length is truncated to $10$ packets.~We begin by describing the policy parametrization that we employ.

\begin{figure}[htbp]
	\begin{centering}
		\includegraphics[trim={1cm 1.6cm 0.9cm 1.3cm},clip, scale = 0.85]{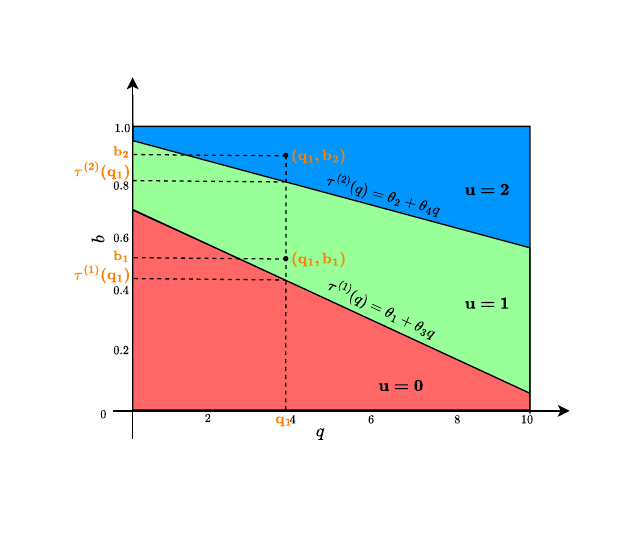} 
		\par\end{centering}
	\caption{Decision regions for $M_d = 2$. For the state $(q,b)$, if we have $b \in [0, \tau^{(1)}(q))$, then it is optimal to transmit $0$ packets, else if $b \in [\tau^{(1)}(q), \tau^{(2)}(q))$, then transmitting $1$ packet is optimal, otherwise if $b \in [\tau^{(2)}(q), 1]$, then it is optimal to transmit $2$ packets. In general, for the state $(q,b)$ if $b \in [\tau^{(j)}(q), \tau^{(j+1)}(q)), j = 0, 1, \ldots, M_d$ with $\tau^0(q) = 0$ and $\tau^{M_d + 1}(q) = 1$, then it is optimal to transmit $j$ packets.}
	\label{fig:regions}
\end{figure}

\begin{figure}[htbp]
	\begin{centering}
		\includegraphics[trim={3.1cm 0.45cm 1.3cm 1.4cm},clip, scale = 0.33]{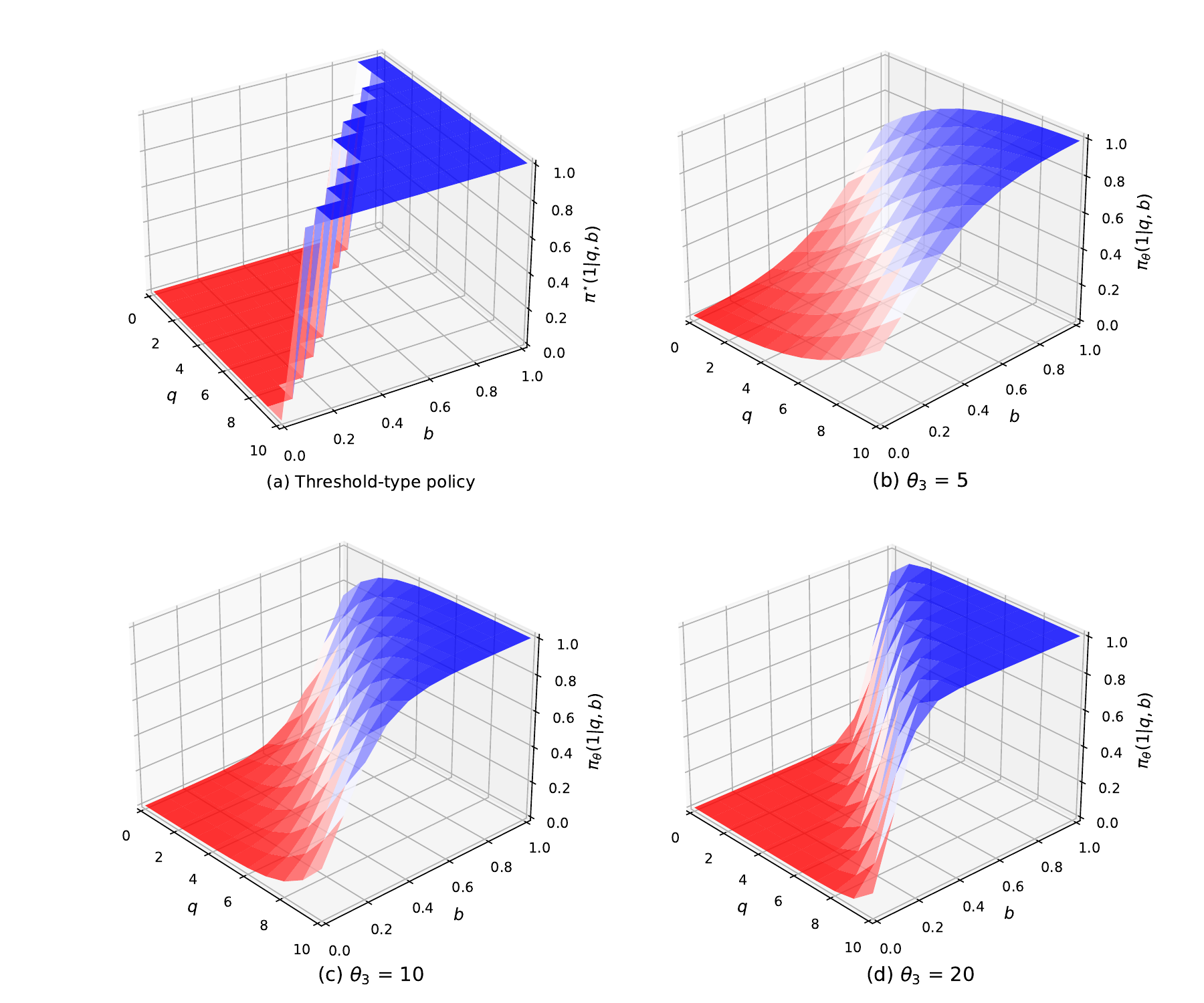} 
		\par\end{centering}
	\caption{(a) Threshold-type policy for $M_d = 1$ with $\tau^{(1)}(q) = 0.9 - 0.08q$, and (b), (c), (d), are its approximation for different values of $\theta_3$ for the parameterized policy $\pi_{\theta}(1|q,b) = \frac{1}{1 + \exp(- {\te_3} (b - \tau^{(1)}(q)))}$.}
	\label{fig:approx_thres_curv}
\end{figure}

\begin{figure}[htbp]
	\begin{centering}
		\includegraphics[trim={5.1cm 0.9cm 3cm 2cm},clip, scale = 0.45]{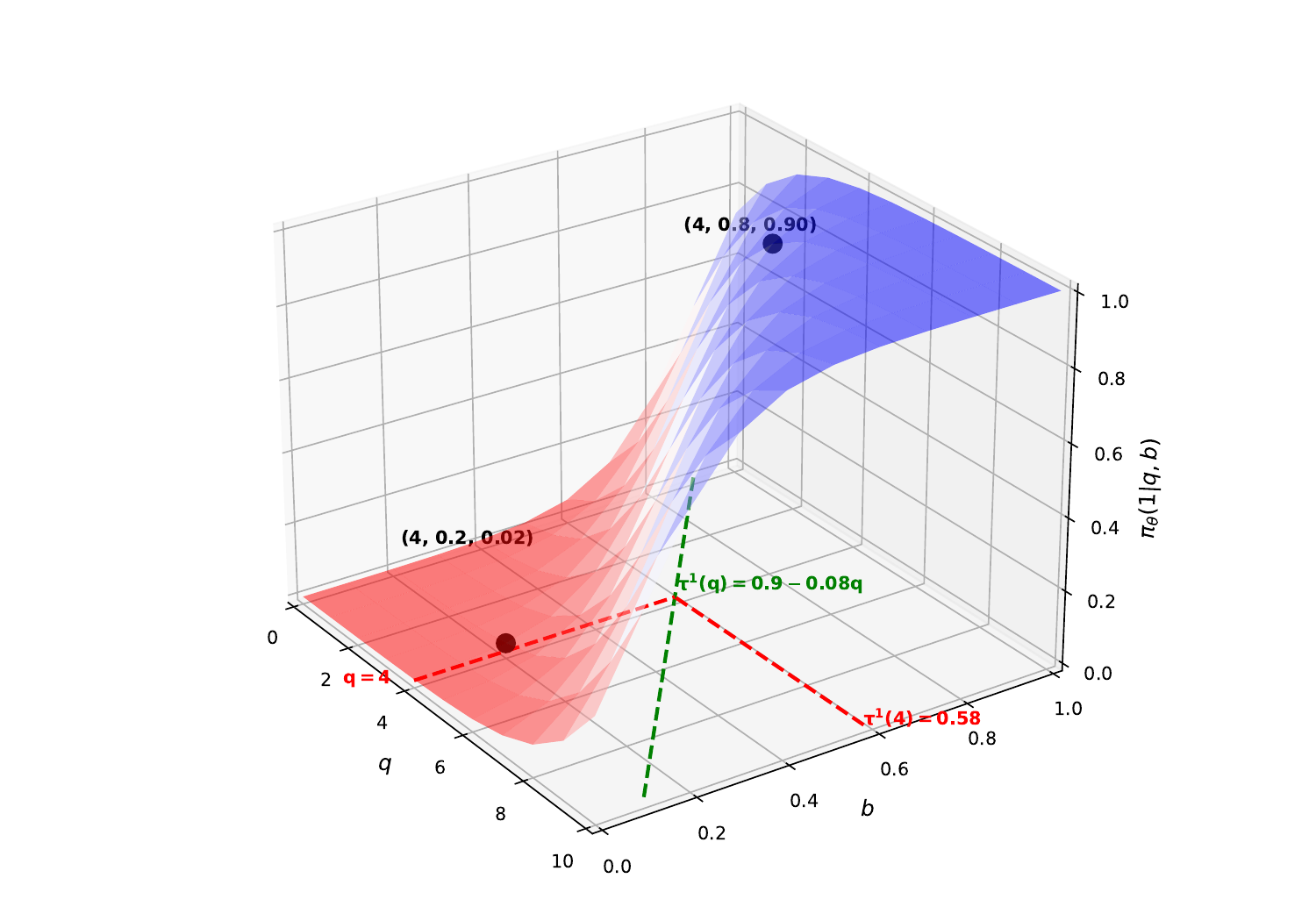} 
		\par\end{centering}
	\caption{Explanation of the parameterized policy for $M_d =1, \tau^{(1)}(q) = \te_1 + \te_2 q$, and $\pi_{\theta}(1|q,b) = \frac{1}{1 + \exp(- {\te_3} (b - \tau^{(1)}(q)))}$ with $\te_1 = 0.9, \te_2 = -0.08$, and $\te_3 = 10$.}
	\label{fig:example}
\end{figure}

\begin{figure}[htbp]
	\begin{centering}
		\includegraphics[scale = 0.29]{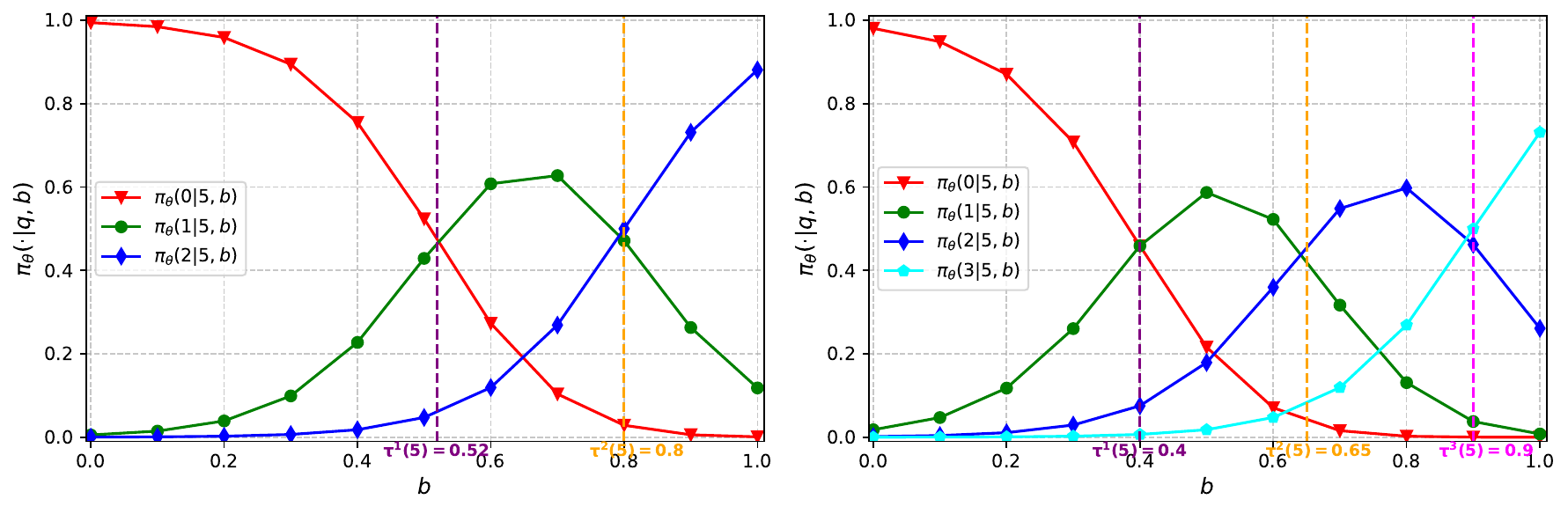} 
		\par\end{centering}
	\caption{Graphical representation of~\eqref{eq:mu_theta_Md-1} for a fixed value of $q = 5$ with (a) $M_d = 2, \tau^{(1)}(q) = 0.92 - 0.08q, \tau^{(2)}(q) = 1 - 0.04q$, and $\theta_5 = \theta_6 = 10$; (b) $M_d = 3, \tau^{(1)}(q) = 0.7 - 0.06q, \tau^{(2)}(q) = 0.9 - 0.05q, \tau^{(3)}(q) = 1 - 0.02q$, and $\theta_7 = \theta_8 = \te_9 = 10$.}
	\label{fig:pi_explain}
\end{figure}

\textbf{Policy parametrization}: We will use the structural result of Theorem~\ref{th:avg_thresh} to perform an efficient search over the policy space. Since a maximum of $M_d$ packets can be transmitted at each time, it follows from Theorem~\ref{th:avg_thresh} that a threshold policy will be characterized by $M_d$ decision regions, or equivalently decision boundaries, where the $j$-th region comprises of those states in which transmitting $j$ packets is optimal. We approximate these decision boundaries by linear curves. More specifically, the $j$-th boundary is given by $\tau^{(j)}(q) = \theta_j + \theta_{M_d + j} q$, where $\theta_j \in \bR, j=1,2,\ldots,2M_d$ are the parameters that completely describe the decision regions. See Fig.~\ref{fig:regions} for an example.
For $j=1,2,\ldots,M_d$ and $(q,b) \in \cX$ define,
\nal{
f^{(\theta)}(q,b;j) :&= \frac{1}{1 + \exp(- (b - (\te_j + \te_{M_d + j} q )  ) \te_{2M_d + j} )}. \notag
}
Let $\pi_{\theta}(j|q, b)$ be the probability with which the policy with parameter $\te= (\te_1,\te_2,\ldots,\te_{3M_d})$ transmits $j$ packets in state $(q,b)$.~Define,
\al{
& \pi_{\theta}(j|q, b) : = f^{(\theta)}(q,b;j) \prod_{i = 1}^{M_d - j} (1 - f^{(\te)}(q,b;j+i)), \label{eq:mu_theta_Md-1}\\
& j = 1,2,\ldots,M_d-1,\notag
}
and let 
\nal{
\pi_{\theta}(M_d|q, b) : = f(q, b; M_d).
}
It is easily verified that the functions $\{\pi_{\te}(j|q,\cdot)\}$ are differentiable, and that $0 \leq \sum_{j = 1}^{M_d} \pi_{\theta}(j|q,b) \leq 1$. We let $\pi_{\theta}(0|q, b) : = 1 - \sum_{j = 1}^{M_d} \pi_{\theta}(j|q,b)$.~From~\eqref{eq:mu_theta_Md-1} we observe that for a state $(q,b)$ that satisfies $b \in [\tau^{(j)}(q),\tau^{(j+1)}(q))$, the term $f^{(j)}(q,b) \rightarrow 1$ as $\theta_{2M_d + j} \rightarrow \infty$; see Fig.~\ref{fig:approx_thres_curv}.~Moreover, the terms $f^{(j+i)}(q,b) \rightarrow 0$ as $\theta_{2M_d + i + j} \rightarrow \infty$.~Thus, $\pi_{\theta}(j|q,b) \rightarrow 1$ as $\te_{2M_d + k}\to \infty$, for all $k =j,j+1,\ldots,3M_d$. This is also demonstrated graphically in Figs.~\ref{fig:example} and~\ref{fig:pi_explain}. We will use the AC algorithm to minimize the average cost over this parametric family of policies $\{\pi_{\theta}, \theta \in \bR^{3M_d}\}$. Let $\zeta(\pi_{\theta})$ be the average cost~\eqref{def:obj} corresponding to policy $\pi_{\te}$.

\textbf{Critic:} Define the advantage function~\cite{sutton1999policy} $\rho_{\theta}$ as follows, with $\tilde{h}_{\theta}$ as defined in~\eqref{eq:4.1.1.a},~\eqref{eq:4.1.1.b}, $(q,b) \in \{0,1,\ldots,10\} \times \cB$, and $u \in \cU$,
\nal{
\rho_{\theta}(q,b,u) = \Tilde{h}_{\theta}(q,b,u) - \zeta(\pi_{\theta}), \forall \theta \in \bR^{3M_d}.
} 
The role of the critic is to estimate the advantage function~\cite{sutton1999policy} associated with the actor's policy. This is then used for tuning the actor's policy parameter $\te$ by updating the actor parameter in an approximate gradient direction of the overall cost.~Critic uses a linear function approximation~\cite{konda1999actor} in order to represent the advantage function. Let $\phi_{\theta}^{(i)}(q,b,u), i = 1,2, \ldots 3M_d$ denote the features used by the critic.~The features $\phi_{\theta}^{(i)}(q,b,u)$ must be chosen such that their span contains the span of the vectors $\{\nabla^{(i)}_{\theta} \log \pi_{\theta}(u|q,b), i = 1, 2, \ldots, 3M_d\}$, where $\nabla^{(i)}_{\theta} \log \pi_{\theta}(u|q,b)$ is the $i$-th component of $\nabla_{\theta} \log \pi_{\theta}(u|q,b)$~\cite{konda1999actor}. One such choice,~\cite{konda1999actor} which we use throughtout is as follows,
\nal{
\phi_{\theta}(q,b,u) = \nabla_{\theta} \log \pi_{\theta}(u|q,b),
}
where $\nabla_{\theta}$ denotes the gradient with respect to $\theta$.~Denote $\phi_{\theta}(q,b,u) = (\phi^{(1)}_{\theta}(q,b,u), \phi^{(2)}_{\theta}(q,b,u), \ldots, \phi^{(3M_d)}_{\theta}(q,b,u))$ $\in \bR^{3M_d}$.~Note that the feature vector $\phi_{\theta}$ depends on the actor parameter vector $\theta$.~Critic updates are performed according to the TD(1) algorithm~\eqref{eq:critic_update}. Then, then the advantage function $\rho_{\theta}$ is approximated as follows,
\al{\label{eq:zeta_approx}
\tilde{\rho}^{(\theta)}_{w} = \sum_{i = 1}^{3M_d} w^{(i)} \phi^{(i)}_{\theta},
}
where $w = (w^{(1)}, w^{(2)}, \ldots, w^{(3M_d)}) \in \bR^{3M_d}$ is the parameter vector for the critic and $\tilde{\rho}$ is an approximation of the true advantage function.

The pseudocode for AC algorithm~\cite{konda1999actor} is presented in Algorithm~\ref{algo:AC}. We use $w(t)$~\eqref{eq:critic_update}, $\theta(t)$~\eqref{eq:actor_update} to denote the critic and actor parameters at time $t$, respectively. Moreover, we use $\hat{R}(t)$~\eqref{eq:zeta} to denote the estimate of the average reward of $\pi_{\theta(t)}$ at time $t$ which is the negative of the average cost, and $z(t)$~\eqref{eq:z} to denote the eligibility trace~\cite{sutton1999policy} at time $t$ in the AC algorithm.
\begin{algorithm}
	\caption{Actor-Critic Algorithm} \label{algo:AC}
	\begin{algorithmic}[0]
            \State Input: policy parameterization $\pi_{\theta}(u|q,b)$
            \State Input: advantage function parameterization $\tilde{\rho}^{(\theta)}_{w}(q,b,u)$
            \State Parameters: actor step size $\alpha^{(\theta)} >0$, \\
            \hspace{1.7cm} critic step size $\alpha^{(w)} > 0$
            \State Initialize $\theta(0)$, $w(0)$, $(Q(0),b(0))$ $ \in $ $\{0,1,\ldots, 10\} \times \cB,  \hat{R}(0)$
            \State Sample $u(0) \sim \pi_{\theta(0)}(\cdot|Q(0),b(0))$ 
            \State Evaluate $z(0) = \phi_{\theta(0)}(Q(0),b(0),u(0))$
		\For {$t =1,2,\ldots, T$}
			\State Take action $u(t-1)$
                \State Observe reward $\bigl[(10 + \kappa c(M_d))- (Q(t-1) + \kappa c(u(t$ 
                \State $-1)))\bigr]$, and next state $(Q(t),b(t))$ as described by~\eqref{eq:2.1.2} 
                \State and~\eqref{eq:2.2.2}
                \State Sample $u(t) \sim \pi_{\theta(t-1)}(\cdot|Q(t),b(t))$
                \State Update average reward: 
                \al{
                & \hat{R}(t) = \hat{R}(t-1) + \alpha^{(w)}\bigl[(10 + \kappa c(M_d)) \notag \\
                &- (Q(t-1) + c(u(t-1))) - \hat{R}(t-1) \bigr] \label{eq:zeta}
                }
                \State Update critic: 
                \al{
                & w(t) = w(t-1) + \alpha^{(w)}\bigl[(10 + \kappa c(M_d))- (Q(t-1) \notag \\
                +& \kappa c(u(t-1))) - \hat{R}(t-1) + \tilde{\rho}^{(\theta(t-1))}_{w(t-1)}(Q(t),b(t),u(t)) \notag \\ 
                -& \tilde{\rho}^{(\theta(t-1))}_{w(t-1)}(Q(t-1),b(t-1),u(t-1))\bigr]z(t-1) \label{eq:critic_update}
                }
                \State Update $z$ for TD(1) critic: Let $(\tilde{q}, \tilde{b}) = (2, p_{11}) \in$
                \State $\{0,1,\ldots, 10\} \times [0,1]$. \al{z(t) = 
                \begin{cases} \label{eq:z}
                    z(t-1) + \phi_{\theta(t-1)}(Q(t),b(t),u(t)) \\
                     \hspace{3.6cm}\mbox{if } (Q(t),b(t)) \neq (\tilde{q}, \tilde{b}), \\
                    \phi_{\theta(t-1)}(Q(t),b(t),u(t)) \hspace{0.15cm}\mbox{otherwise} 
                \end{cases}}
                \State Update actor: 
                \al{
                \theta(t) &= \theta(t-1) \notag\\
                &+ \alpha^{(\theta)} \tilde{\rho}^{(\theta(t-1))}_{w(t-1)}(Q(t),b(t),u(t)) z(t) \label{eq:actor_update}
                }
		\EndFor
	\end{algorithmic} 
\end{algorithm}
As is shown in~\cite[Theorem 2]{konda1999actor}, under Assumptions (A1)-(A6) of~\cite{konda1999actor}, the AC algorithm is guaranteed to converge to a local minima of the average cost function $\zeta(\pi_{\cdot})$, so that this gives us a locally optimal policy for the POMDP~\eqref{def:pomdp}. This result is stated next.

\begin{theorem}\label{thm:AC_convg}
    Consider the AC algorithm with TD(1) critic (Algorithm~\ref{algo:AC}) used for solving the POMDP~\eqref{def:pomdp}, in which the actor tunes the policy parameter $\te(t)$, while the critic tunes the weights $w(t)$ that are used for approximating the advantage function according to~\eqref{eq:zeta_approx}. Then, we have that $\liminf_{t} ||\nabla_{\theta} \zeta(\pi_{\theta(t)})|| = 0$ with probability $1$, where $\zeta(\cdot)$ is the average cost function and $\nabla_{\theta} \zeta(\pi_{\theta(t)})$ is the gradient of $\zeta(\pi_{\theta})$ with respect to $\theta$ evaluated at $\theta(t)$.
\end{theorem}

\emph{Simulation setup:} Throughout, we set the initial state of the POMDP to be $(Q(0),b(0)) = (5,0.5)$. Markovian channel's parameters are set equal to $p_{01} = 0.4,~p_{11} = 0.9$, and packet transmission cost function is taken to be $c(u) = \exp(u)-1$. For the AC algorithm, we initialize the estimate of average reward as $\hat{R}(0) = 0$~\eqref{eq:zeta}. Initial value of the Actor's parameter $\theta(0)$
and critic's parameter $w (0)$~\eqref{eq:critic_update} are chosen uniformly at random from $[0,1]^{3M_d}$. 

\begin{figure}[htbp]
	\begin{centering}
		\includegraphics[scale=0.35]{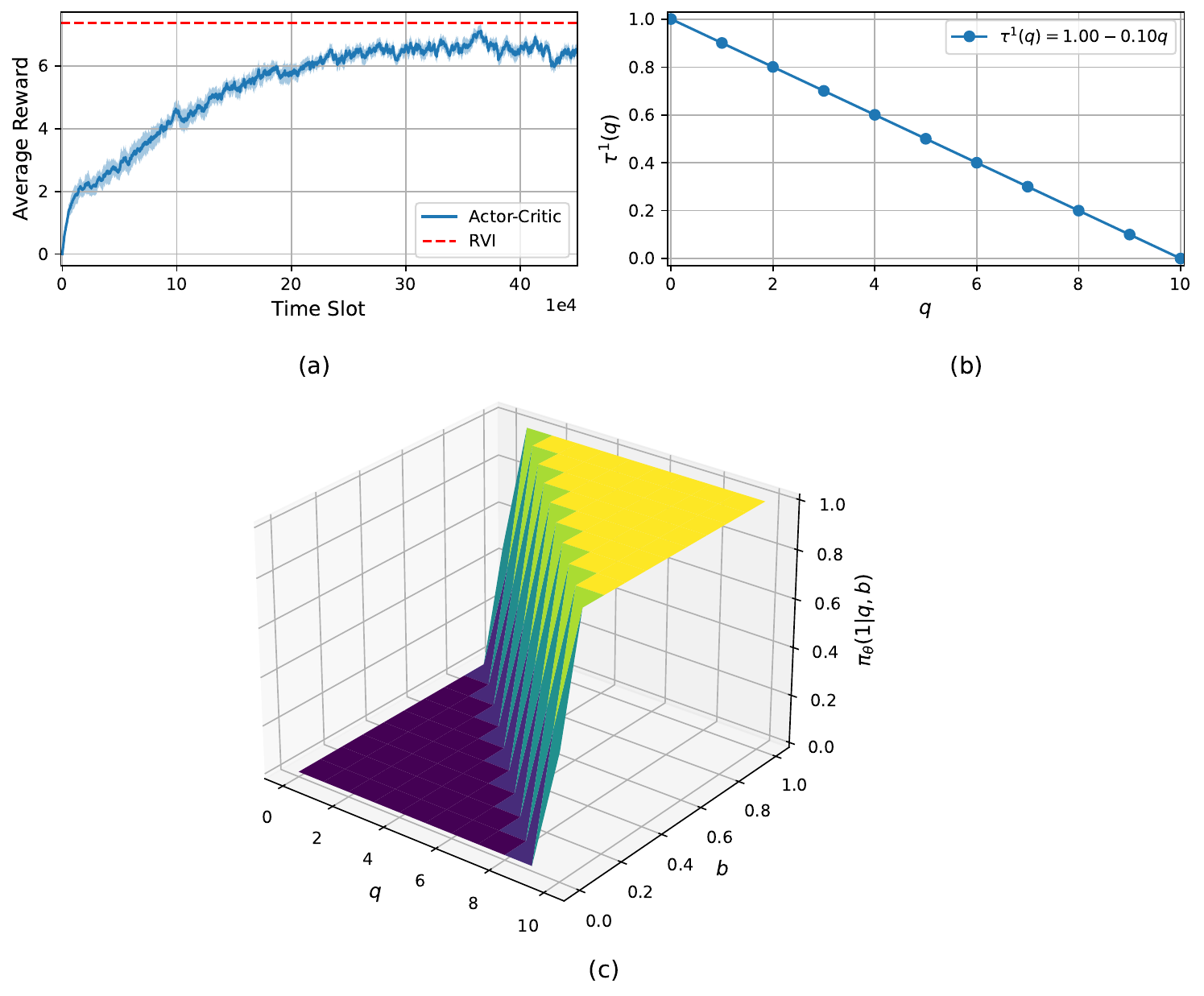} 
		\par\end{centering}
	\caption{(a) Comparison of the average reward achieved by the AC and RVI with $T=450000,~\alpha^{(\theta)} = 0.0005$, $\alpha^{(w)} = 0.002,~p_{1} = 0.7,~p_{01} = 0.4$, and $p_{11} = 0.9$: (b) Threshold curve, $\tau^{(1)}(q) = 1 - 0.1q$ obtained using the AC; (c) Policy obtained from the AC.}
	\label{fig:avg_rwd}
\end{figure}

\begin{figure}[htbp]
	\begin{centering}
		\includegraphics[scale=0.35]{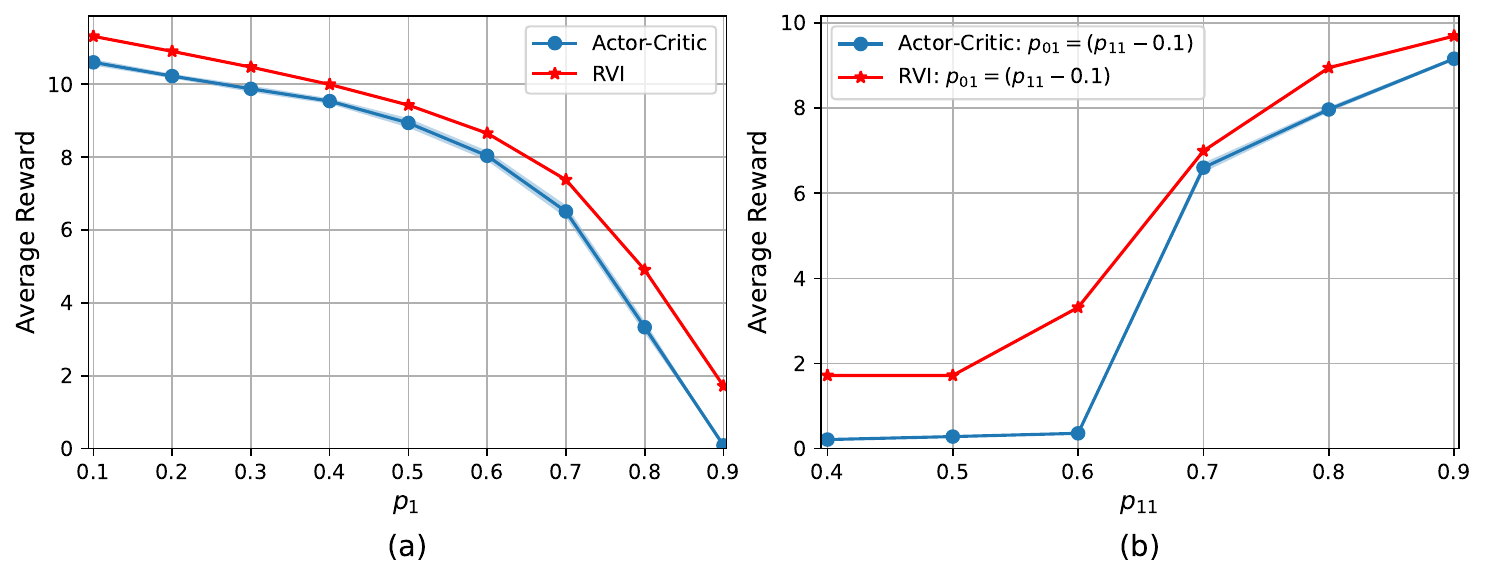} 
		\par\end{centering}
	\caption{Performance comparison of the AC and RVI as network parameters are varied. We use $T=450000$, while step-sizes of the AC are held fixed at $\alpha^{(\theta)} = 0.0005$, and $\alpha^{(w)} = 0.002$: (a) $p_{01} = 0.4, p_{11} = 0.9$; (b) $p_{1} = 0.6$.}
	\label{fig:change_in_params}
\end{figure}

\textbf{Case (i): $M_a = 1$ and $M_d =1$}: In this case we have just one decision boundary, so we use three parameters $\te_1,\te_2,\te_3$ to parameterize this.~Fig.~\ref{fig:avg_rwd}(a) compares the performance of the AC with the RVI.~In fact, it can be observed that the terminal value of average reward is almost equal to the optimal average reward that is obtained by the RVI algorithm.~Fig.~\ref{fig:avg_rwd}(b) depicts the threshold curve $\tau^{(1)}(q) = 1 - 0.1q$ attained by the AC. $\tau^{(1)}(\cdot)$ is monotone decreasing.~Note that we have not formally proved that the threshold curve $\tau^{(1)}(\cdot)$ will be monotone, this could be the subject of a future work.~Fig.~\ref{fig:avg_rwd}(c) shows the policy obtained from the AC. Specifically, it plots the probability with which a packet is transmitted by the transmitter.~Next, in Fig.~\ref{fig:change_in_params} we compare the performance of the AC and RVI as the network parameters are varied. It can be seen from Fig.~\ref{fig:change_in_params}(a) that, as the packet arrival probability $p_{1}$ increases, the average reward decreases. This is because as $p_{1}$ increases, the queue length increases, leading to higher instantaneous costs and more frequent transmissions to manage the queue. Both of these result in a decrease in the average reward. Moreover,~Fig.~\ref{fig:change_in_params}(b) indicates that as the channel becomes better, i.e. $p_{11}$ increases, the average reward increases. Moreover, the optimal policy is seen to outperform the learning policy.


\begin{figure}[htbp]
	\begin{centering}
		\includegraphics[scale=0.35]{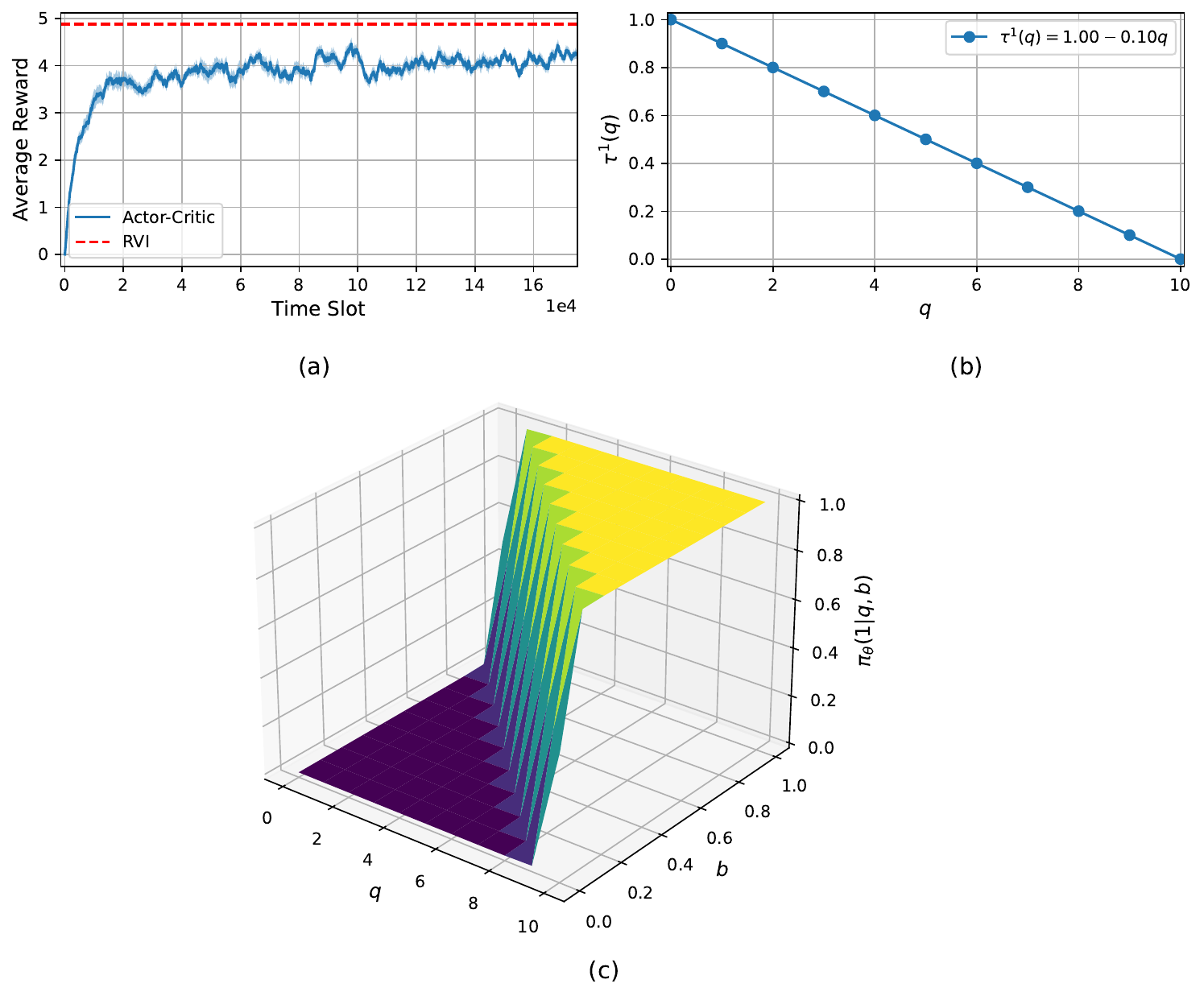} 
		\par\end{centering}
	\caption{(a) Comparison of the average reward achieved by the AC and RVI with $T= 175000,~\alpha^{(\theta)} = 0.0003$,~$\alpha^{(w)} = 0.002$, $M_a = 2,~M_d = 1,~p_{01} = 0.4,~p_{11} = 0.9$, $p_{1} = 0.4$, and $p_{2} = 0.2$ ; (b) Threshold curve attained by the AC, $\tau^{(1)}(q) = 1 - 0.01q$; (c) Policy obtained from AC.}
	\label{fig:avg_rwd_2_pckts}
\end{figure}

\textbf{Case (ii): $M_a = 2$ and $M_d = 1$}: Fig.~\ref{fig:avg_rwd_2_pckts}(a) compares the performance of the AC and   RVI.~Fig.~\ref{fig:avg_rwd_2_pckts}(b) shows that $\tau^{(1)}(q)$ is a monotone decreasing function of the queue length $q$. Fig.~\ref{fig:avg_rwd_2_pckts}(c) shows the policy obtained from the AC.~We see that even though the AC is guaranteed to converge only to a local maxima, its performance is almost the same as RVI.~Also, it suffices to use simple linear decision boundaries in order to achieve a near-optimal performance given by RVI.

\begin{figure}[htbp]
	\begin{centering}
		\includegraphics[scale=0.35]{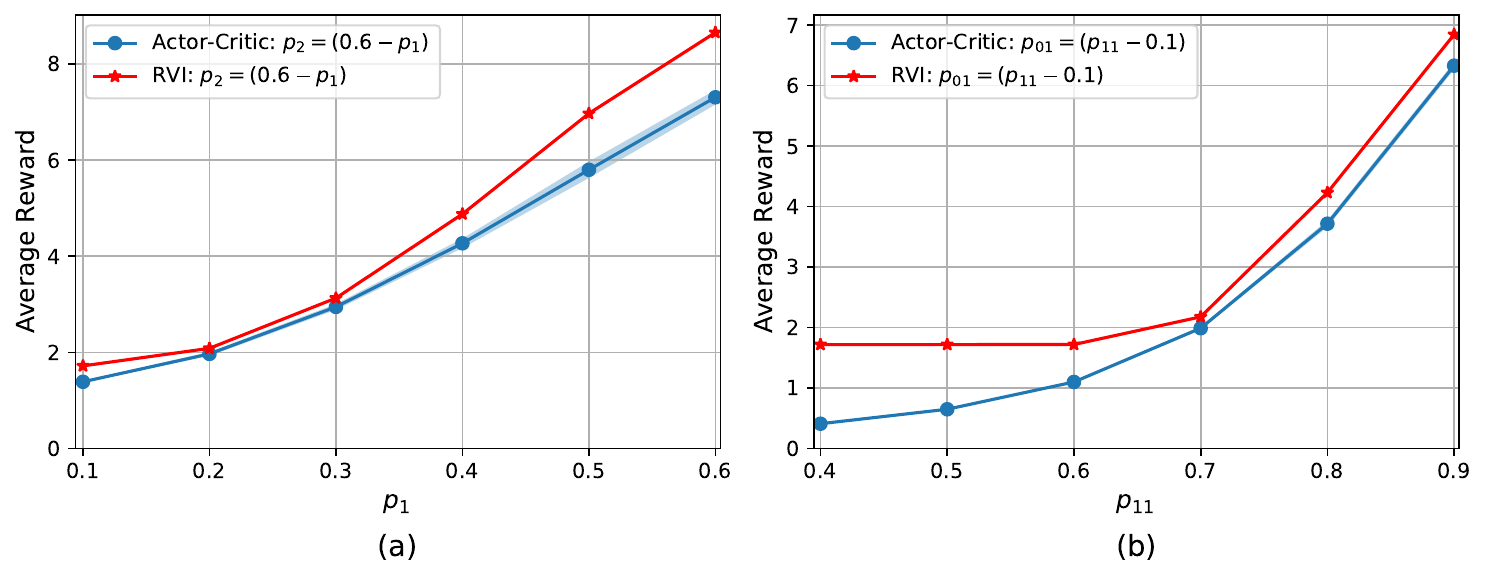} 
		\par\end{centering}
	\caption{Performance comparison of the AC and RVI as network parameters are varied $T=175000,~\alpha^{(\theta)} = 0.0003$,~$\alpha^{(w)} = 0.002$, $M_a = 2$, $M_d = 1$: (a) $p_{01} = 0.4, p_{11} = 0.9$; (b) $p_{1} = 0.4, p_{2} = 0.2$.}
	\label{fig:avg_rwd_params_2_pckts}
\end{figure}

Fig.~\ref{fig:avg_rwd_params_2_pckts} compares the performance of the AC and RVI as the packet arrival probabilities $p_{1}, p_{2}$ are varied. We observe that since the mean packet arrival rate decreases with a decrease in $p_{2}$, the average reward of policy returned by the AC, as well as that of the RVI algorithm, increases.~It can be seen that the AC algorithm yields a performance that is near-optimal since it is almost equal to that of the RVI.

\begin{figure}[htbp]
	\begin{centering}
		\includegraphics[scale=0.34]{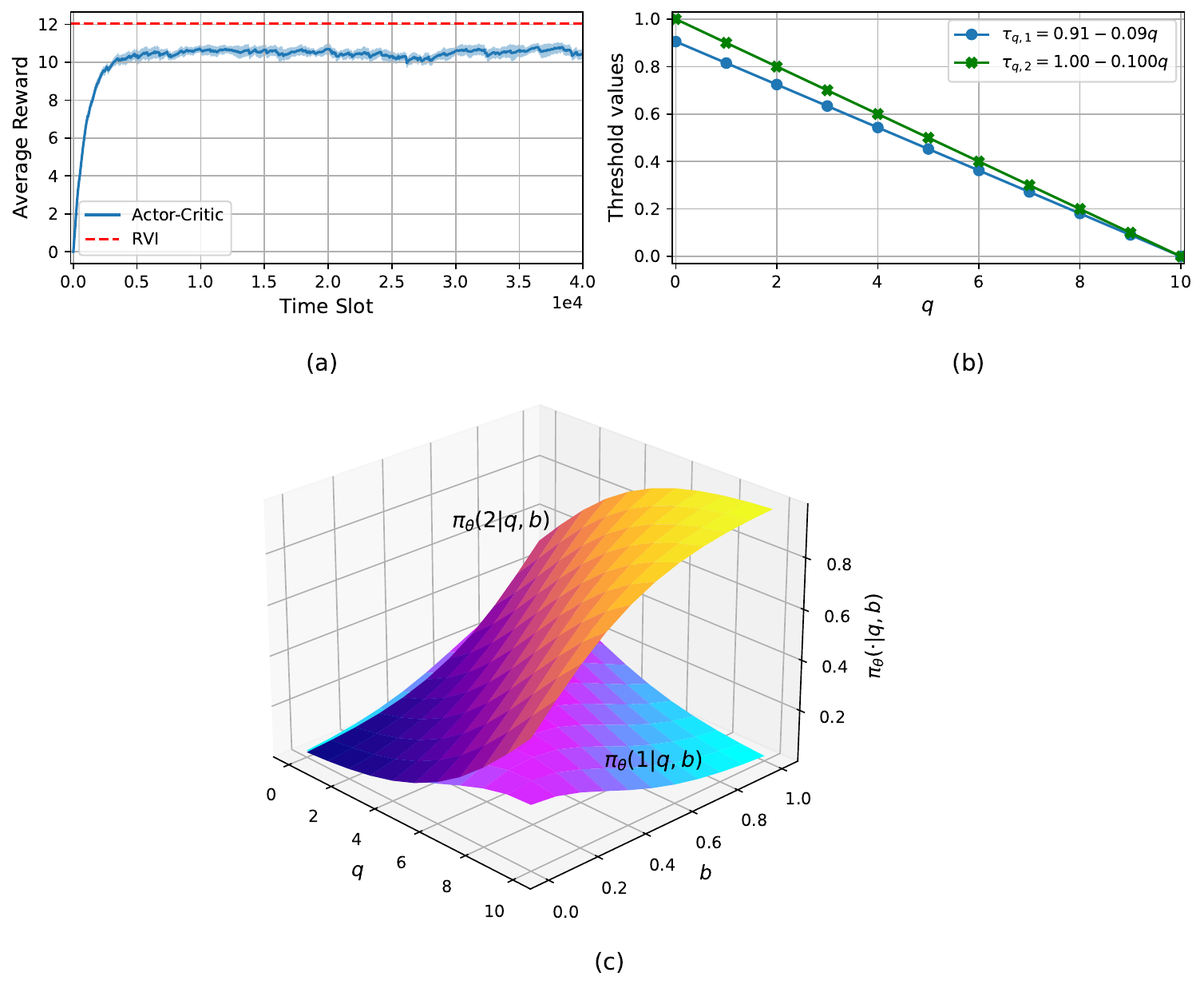} 
		\par\end{centering}
	\caption{(a) Comparison of the average reward achieved by the AC and RVI with $T=40000,~\alpha^{(\theta)} = 0.0006$, $\alpha^{(w)} = 0.001$, $M_a = 1, M_d = 2, p_{01} = 0.4, p_{11} = 0.9$, and $p_{1} = 0.9$; (b) Threshold curves attained by the AC; (c) Policy obtained from the AC.}
	\label{fig:avg_rwd_u_3}
\end{figure}

\begin{figure}[htbp]
	\begin{centering}
		\includegraphics[scale=0.35]{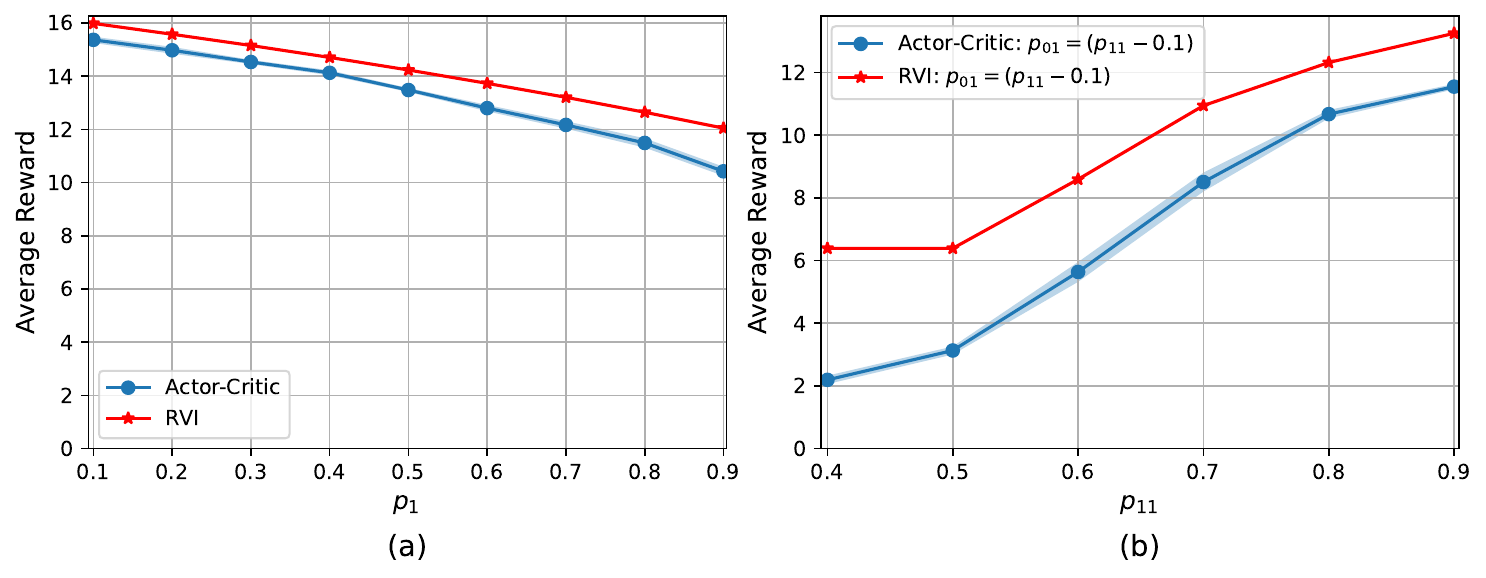} 
		\par\end{centering}
	\caption{Performance comparison of the AC and RVI as network parameters are varied with $T=40000,~\alpha^{(\theta)} = 0.0006$, $\alpha^{(w)} = 0.001$, $M_a = 1, M_d = 2$: (a) $p_{01} = 0.4, p_{11} = 0.9$; (b) $p_{1} = 0.9$.}
	\label{fig:avg_rwd_params_u_3}
\end{figure}

\textbf{Case (iii): $M_a = 1$ and $M_d =2$}: Here $\theta = (\theta_1, \theta_2, \ldots, \theta_6) \in \bR^6$ and there will be three decision regions (two decision boundaries).~Fig. \ref{fig:avg_rwd_u_3}(a) compares the performance of the AC and RVI, while~Fig. \ref{fig:avg_rwd_u_3}(b) depicts the threshold curve of the policy returned by the AC algorithm. It can be seen that $\tau^{(2)}(q) > \tau^{(1)}(q)$, so that more packets are transmitted when the channel is more likely to be in a good condition. This property is shown in Proposition~\ref{prop:2}. Moreover,~Fig. \ref{fig:avg_rwd_u_3}(c) depicts that as the queue length increases, the probability of transmitting $1$ packet decreases, while the probability of transmitting $2$ packets increases. Fig.~\eqref{fig:avg_rwd_params_u_3} compares the performance of the AC and RVI as the network parameters are varied.

\section{Conclusion}\label{sec:conclusion}
We design dynamic scheduling policies for a mmWave network in which a single transmitter sends data packets over a Gilbert Elliott channel.~These minimize an average cost consisting of two components: the average end-to-end packet delay, and the packet transmission cost. The transmitter gets to observe the channel-state only when it attempts a transmission, so that this gives rise to a POMDP.~We show that this POMDP has an optimal policy with a threshold structure. Since threshold policies are simple to implement, this result leads to a practical scheduling solution.~Moreover, we leverage this structural insight to efficiently search within a class of parametric threshold policies using the actor-critic learning algorithm, which is particularly valuable when the system parameters are unknown.

\bibliographystyle{IEEEtran}
\bibliography{mmWaveReference_v2}
\appendices

\section{Existence of an Optimal Stationary Policy}\label{app:stat_optimal}
\begin{theorem}
	There exists a stationary deterministic policy that is optimal for
	the discounted cost POMDP (\ref{def:discountedpomdp}).\label{thm:countable} 
\end{theorem}
\begin{IEEEproof}
	We begin by showing that the state space $\mathcal{X}=\bZ_{+} \times \cB$ associated
	with the POMDP (\ref{def:discountedpomdp}) is countable. The set $\bZ_{+}$ is countable. We will show that the set $\mathcal{B} = \mathcal{B}_{0}\cup\mathcal{B}_{1}\cup\mathcal{B}_{2}$ which contains the belief state, is also countable. To show this, map the $k$-th element
	of $\mathcal{B}_{j}$ to the integer $3k+j$, $j \in \{0,1,2\}$. This is a one-to-one map from $\mathcal{B}_{0}\cup\mathcal{B}_{1}\cup\mathcal{B}_{2}$
	to the set of natural numbers, and hence $\mathcal{B}_{0}\cup\mathcal{B}_{1}\cup\mathcal{B}_{2}$
	is countable. Since $\bZ_{+}$ and $\cB$ are countable, $\bZ_{+} \times \cB$ is also countable. Since for this POMDP the state-space $\mathcal{X}=\bZ_{+} \times \cB$ is countable, and the action space is finite, it follows from~\cite[Th. 6.2.10]{puterman2014Markov}
	that there exists a stationary deterministic policy that is optimal
	for (\ref{def:discountedpomdp}). 
\end{IEEEproof}

\section{Auxiliary Results used in Section~\ref{subsec:General-threshold}}

\begin{lemma}
	The functions $\tilde{J}_{n}\ubeta\left(q,\cdot;0\right)$~\eqref{eq:3.2.1},~\eqref{eq:3.3.1}, $q\in\bZ_{+},~n\in\bN$, are concave.
	\label{lem:lemma2} 
\end{lemma}
\begin{IEEEproof}
	\noindent We prove this by using induction on $n$. The statement holds trivially for $n=1$ since $\tilde{J}_{1}\ubeta\left(q,b;0\right)=q$.~Next, assume that the functions $\tilde{J}_{n}\ubeta\left(q,\cdot;0\right),~q\in \bZ_{+},~n=1,2,\ldots,m$ are concave.~We will show that $\tilde{J}_{m+1}\ubeta\left(q,\cdot;0\right),~q\in \bZ_{+}$ are concave. It follows from (\ref{eq:3.2.1})
	that it suffices to prove that $J_{m}\ubeta\left(q,\mathcal{T}\left(b\right)\right),~q\in \bZ_{+}$
	are concave in $b$. In the remaining proof, we focus exclusively on proving concavity of $J_{m}\ubeta\left(q,\mathcal{T}\left(b\right)\right)$ with respect to $b$ for all $q\in\bZ_{+}$.~It follows from (\ref{eq:3.1})
	and (\ref{eq:3.0.2}) that 
	\begin{equation}
		J_{m}\ubeta\left(q,\mathcal{T}\left(b\right)\right)=\underset{u\in\cU }{\min}~\tilde{J}_{m}\ubeta\left(q,\mathcal{T}\left(b\right);u\right).\notag
	\end{equation}
	Now, $\tilde{J}_{m}\ubeta\left(q,\mathcal{T}\left(\cdot\right);0\right)$
	is concave since it is the composition of an affine function $\mathcal{T}\left(\cdot\right)$,
	and a concave function $\tilde{J}_{m}\ubeta\left(q,\cdot;0\right)$.~For $u\geq1$, upon substituting $\mathcal{T}\left(b\right)=b p_{11}+\left(1-b\right)p_{01}$ in (\ref{eq:3.3.1}), we see that the functions $\tilde{J}_{m}\ubeta\left(q,\mathcal{T}\left(\cdot\right);u\right)$ are affine, and hence concave.~Since $J_{m}\ubeta\left(q,\mathcal{T}\left(b\right)\right)$ is the pointwise
	minima of $\{J_{m}\ubeta\left(q,\mathcal{T}\left(b\right);u\right)\}_{u=0}^{M_d}$, each of which are concave, $J_{m}\ubeta\left(q,\mathcal{T}\left(b\right)\right)$ is also concave. This completes the proof. 
\end{IEEEproof}

\begin{lemma}
	The functions $J_{n}\ubeta\left(q,\cdot\right)$~\eqref{eq:3.0.2}, $q \in \bZ_+,~n\in\bN$ are concave, and hence for each $q\in\bZ_{+}$, $J\ubeta\left(q,\cdot\right)$~\eqref{eq:3.1} is also concave. \label{lem:lemma3} 
    
\end{lemma}
\begin{IEEEproof}
	\noindent Fix a $q \in\bZ_{+}$. From (\ref{eq:3.0.2}) we have that $J_{n}\ubeta\left(\cdot,\cdot\right)$
	is a pointwise minima of the functions $\left\{\tilde{J}_{n}\ubeta\left(\cdot,\cdot;u\right)\right\}_{u\in\cU}$.~It follows from Lemma \ref{lem:lemma2} that $\tilde{J}_{n}\ubeta\left(q,\cdot;0\right)$ is concave. Also, when $u\ge 1$, it follows from (\ref{eq:3.3.1}) that $\tilde{J}_{n}\ubeta\left(q,\cdot;u\right)$ is affine, and hence concave. Since $J_{n}\ubeta\left(q,\cdot\right)$
	is a pointwise minima of concave functions $\left\{\tilde{J}_{n}\ubeta\left(q,\cdot;u\right)\right\}_{u\in\cU}$, it is concave. It then follows from~\eqref{eq:Jn_convg_J} that $J\ubeta\left(q,\cdot\right)$ is also concave. This completes the proof.
\end{IEEEproof}
\begin{lemma}
	The functions $J_{n}\ubeta\left(\cdot,b\right)$~\eqref{eq:3.0.2}, $n\in\bN,b\in\cB$ are monotonic non-decreasing, and hence for each $b\in\mathcal{B}$, $J\ubeta\left(\cdot,b\right)$~\eqref{eq:3.1} is also monotonic non-decreasing.\label{lem:lemma4-1}
	 
\end{lemma}
\begin{IEEEproof}
	We prove this using induction on $n$. The claim holds trivially for $n=0$ since 
	$J_{0}\ubeta\left(q,b\right)=0$.~Next, assume that the functions $J_{n}\ubeta\left(\cdot,b\right),~b\in \mathcal{B}$ are monotonic non-decreasing for $n=0,1,\ldots,m$.~We will show that the function $J_{m+1}\ubeta$ satisfies the following: whenever $q_1>q_2$, we have $J_{m+1}\ubeta\left(q_{1},b\right)\geq J_{m+1}\ubeta\left(q_{2},b\right)$. From (\ref{eq:3.0.2}) we have that $J\ubeta_{m+1}(\cdot,b)$ is a pointwise minima of the functions $\tilde{J}\ubeta_{m+1}(\cdot,b;0), \{\tilde{J}\ubeta_{m+1}(\cdot,b;u)\}_{u=1}^{M_d}$ it suffices to show that the following hold: $\tilde{J}_{m+1}\ubeta\left(q_{1},b;0\right)\geq \tilde{J}_{m+1}\ubeta\left(q_{2},b;0\right)$
	and $\tilde{J}_{m+1}\ubeta\left(q_{1},b;u\right)\geq \tilde{J}_{m+1}\ubeta\left(q_{2},b;u\right), u\in \{1,2,\ldots,M_{d}\}$.
	We have,
	\nal{
	\tilde{J}_{m+1}\ubeta\left(q_{1},b;0\right) & =q_{1}+\beta\stackrel[i=0]{M_{a}}{\sum}p_{i}J_{m}\ubeta\left(\cQ^+(q_1,i,0),\mathcal{T}\left(b\right)\right) \\
	& \geq q_{2}+\beta\stackrel[i=0]{M_{a}}{\sum}p_{i}J_{m}\ubeta\left(\cQ^+(q_2,i,0),\mathcal{T}\left(b\right)\right) \\
	& =\tilde{J}_{m+1}\ubeta\left(q_{2},b;0\right),
	}
	where the inequality is true since from the induction hypothesis we have that $J_{m}\ubeta\left(\cdot,b\right)$
	is monotonic non-decreasing.~Similarly, for $u\in \{1,2,\ldots,M_{d}\}$, we have 
	
	\begin{eqnarray}
		&  & \tilde{J}_{m+1}\ubeta\left(q_{1},b;u\right)\nonumber \\
		& = & q_{1}+ \kappa c(u)+\beta\stackrel[i=0]{M_{a}}{\sum}p_{i}\left(b J_{m}\ubeta\left(\cQ^+(q_1,i,u),p_{11}\right)\right.\nonumber \\
		&  & \left.\vphantom{\beta\left(b J_{n}\ubeta\left(\left(q+i-j\right)_{0}^{M_{q}},p_{11}\right)+\right.}+\left(1-b\right)J_{m}\ubeta\left(\cQ^+(q_1,i,0),p_{01}\right)\right)\nonumber 
        \end{eqnarray}
    \begin{eqnarray}
		& \geq & q_{2}+ \kappa c(u)+\beta\stackrel[i=0]{M_{a}}{\sum}p_{i}\left(b J_{m}\ubeta\left(\cQ^+(q_2,i,u),p_{11}\right)\right.\nonumber \\
		&  & \left.\vphantom{\beta\left(b J_{n}\ubeta\left(\left(q+i-j\right)_{0}^{M_{q}},p_{11}\right)+\right.}+\left(1-b\right)J_{m}\ubeta\left(\cQ^+(q_2,i,0),p_{01}\right)\right) \notag \\
		& = & \tilde{J}_{m+1}\ubeta\left(q_{2},b;u\right),\nonumber 
	\end{eqnarray}
	where the inequality holds because $J_{m}\ubeta\left(q,b\right)$
	is monotonic non-decreasing in $q$ based on the induction hypothesis. It then follows from~\eqref{eq:Jn_convg_J} that $J\ubeta\left(\cdot,b\right)$ is also monotonic non-decreasing. This completes the proof. 
\end{IEEEproof}
\section{Auxiliary Results used in Section \ref{subsec:Optimal-Policy-Average}} \label{app:Boundedness-Proof}
Define $\widetilde{\mathcal{X}} \subset \mathcal{X}$ as follows,
\begin{equation}
	\widetilde{\mathcal{X}}:=\underset{q\in\mathbb{Z}^{+}}{\bigcup}\left\{ \left(q,p_{01}\right),\left(q,p_{11}\right)\right\} .\notag
\end{equation}
We begin by showing some properties of the Markov process $\{(Q(t),b(t))\}$.
\begin{lemma}\label{lem:irreducible}
Consider the system described by~\eqref{eq:2.1.2},~\eqref{eq:2.2.2}, and the Markov process $\{(Q(t),b(t))\}$ which results from transmitting $\min\{Q(t),M_{d}\}$ packets at each time $t$.~The state-space of this Markov chain is $\widetilde{\mathcal{X}}$. Under Assumption \ref{assumption2}, this Markov process is irreducible and aperiodic. 
\end{lemma}
\begin{IEEEproof}
We begin by proving irreducibility. We will show that given any two states $x,y \in \widetilde{\mathcal{X}}$, there exists a path connecting $x$ to $y$ and that has a strictly positive probability. Fix a sufficiently large time duration $T_0$.~$p_{0}>0$ is the probability that no packet arrives. Consider the event that no packet arrival occurs for $T_0$ steps. Moreover, since the Markovian channel probabilities $p_{01},p_{11}>0$, and $T_0$ is sufficiently large, with a positive probability the cumulative number of packet deliveries during $T_0$ is equal to the initial queue length. Hence, the state value after $T_0$ is $\left(0,p_{11}\right)$. 

Next, we will show that beginning from $\left(0,p_{11}\right)$, it is possible to reach any state in $\widetilde{\mathcal{X}}$ with a positive probability.~In order to reach $(q,p_{01})$, the process firstly hits the state $(0,p_{01})$ (this occurs with positive probability since $p_{0}>0$
	and $p_{11}<1$), and then it reaches $\left(q,p_{01}\right)$ (this happens with a positive probability since $p_{1}>0$ and
	$p_{01}<1$). Finally, in order to reach $\left(q,p_{11}\right)$ from $(0,p_{11})$, the process firstly hits $\left(q+M_{d},p_{01}\right)$ in finite time (as discussed above) and then $\left(q,p_{11}\right)$ in finite time (since $p_{0}>0$ and $p_{01}>0$).~This shows that the Markov chain is irreducible.~Aperiodicity follows since any state of the form $\left(q,p_{01}\right)$ has a self loop since $p_{0}>0$ and $p_{01}<1$. This completes the proof.
\end{IEEEproof}
Recall that $\mu = (\mu_0,\mu_1)$ is the vector consisting of stationary probabilities associated with the channel state process $\{s(t)\}$.
\begin{lemma}\label{lem:average-cost-finite}
Under Assumption~\ref{assum:stability}, the average cost~\eqref{def:pomdp} incurred by the policy which transmits $\min\{Q(t),M_{d}\}$ packets at each time $t$, is finite.
\end{lemma}
\begin{IEEEproof}
Consider a time horizon $T_0 \in \bN$ that is much greater than the mixing time of the Markov process associated with the channel state process $\{s(t)\}$. Consider the Markov chain $\{Q(t)\}_{t\in\bN}$ which results by the application of the policy discussed above. We will show that this chain is positive recurrent. It is shown in Lemma~\ref{lem:irreducible} that $\{Q(t)\}_{t\in\bN}$ is irreducible and aperiodic.~Consider the Lyapunov function $Q(t)^2$.~The Lyapunov drift~\cite{meyn2012markov} under this policy satisfies the following bound
	\begin{align}
		&\bE\left\{ Q(kT_0)^2 - Q((k-1)T_0)^2| \cF_{(k-1)T_0} \right\} \le     T^2_0\notag\\
		&+ 2 Q(kT_0) \bE\left(\sum_{t=(k-1)T_0}^{kT_0} A(t) - D(t)     \Bigg| \cF_{(k-1)T_0}\right),  \label{ineq:drift}
	\end{align}
where $D(t)$ is the number of packets that are successfully delivered at $t$, while $\cF_{\ell}$ is the sigma-algebra generated by the random variables $\{u(k)\}_{k=1}^{l-1},\{A(k)\}_{k=1}^{l},\{u(k)s(k)\}_{k=1}^{l-1}$. Now, if $Q((k-1)T_0)> T_0 M_d$, then the queue can never be empty in the steps $(k-1)T_0,(k-1)T_0+1,\ldots,kT_0$, so that $D(t)$ is equal to $M_d$ if the channel state $s(t)$ is equal to $1$, and is $0$ otherwise. Since $T_0$ is much larger than the mixing time of $\{s(t)\}$, this yields us,    
    \nal{
    & \bE \left(\sum_{t=(k-1)T_0}^{kT_0} D(t) \Bigg| \cF_{(k-1)T_0}\right) \\
    & = M_d \sum_{j=0}^{T_0} \mathbb{P}\left(s((k-1)T_0 + j) = 1 | \cF_{(k-1)T_0}\right) \\
    & = M_d \sum_{j=0}^{T_0} \left(\mu_1 - (p_{11} - p_{01})^j (\mu_1 - b_0) \right) \\
    & = M_d \left(\mu_1 T_0 - (\mu_1 - b_0) \frac{1- (p_{11} - p_{01})^{T_0 + 1}}{p_{10} + p_{01}}\right) \\
    & = M_d\left(\mu_1 T_0 - \frac{\mu_1 - b_0}{p_{10} + p_{01}}\right) + \frac{\mu_1 - b_0}{p_{10} + p_{01}} \\
    & \times \exp\{-(-\ln{(p_{11}-p_{01})})\} \exp\left\{-(-\ln{(p_{11}-p_{01})})) T_0\right\},
    }
where $b_0 = \mathbb{P}\left(s((k-1)T_0) = 1 | \cF_{(k-1)T_0}\right)$ and the second equality follows from~\cite[Lemma 1]{pomdp_whittle1}. Moreover, since arrivals are i.i.d., we have $\bE \sum_{t=(k-1)T_0}^{kT_0} A(t) = T_0 \bE A(1)$. Substituting these into~\eqref{ineq:drift} we obtain
	\begin{align}
		&\bE\left\{ Q(kT_0)^2 - Q((k-1)T_0)^2| \cF_{(k-1)T_0} \right\} \notag \\ & \le   
        T_0 \left(T_0 + 2 Q(kT_0) \left[\bE A(1) - M_d \mu_1 \right]   \right)
,\label{ineq:drift1}
\end{align}
so that $\bE\left\{ Q(kT_0)^2 - Q((k-1)T_0)^2| \cF_{(k-1)T_0} \right\} <0$ for sufficiently large values of $Q(kT_0)$ since we had assumed $M_d \mu_1 > \bE A(1)$.~From Foster Lyapunov Theorem~\cite{meyn2012markov} we have that $\{Q(t)\}_{t\in\bN}$ is positive recurrent, and admits a unique stationary distribution. Now, consider the process $\{Q(t)\}_{t\in\bN}$ starting with an initial distribution equal to the stationary distribution, and let $Q(\infty)$ be a random variable that has this distribution. We would then have that $\bE Q(kT_0)^2 = \bE Q((k-1)T_0)^2$, so that upon
	using law of total expectation in~\eqref{ineq:drift1}, and substituting this into the resulting expression, we obtain the following upper-bound on mean value of queue length,
	\begin{align}
		\bE Q(\infty)  \le     \frac{T_0}{2\left(M_d \mu_1 - \bE A(1) \right)}.\notag
	\end{align}
This completes the proof since the transmission cost is bounded.
\end{IEEEproof}

The following lemmas validate assumptions $1$, $2$, $3$, and $3^{\star}$ of \cite{sennott1989average} and thereby help to extend the threshold structure of the discounted
cost POMDP (\ref{def:discountedpomdp}) to the average cost~(\ref{def:pomdp}).
\begin{lemma}
	\label{lem:SennottC1}For every state $\left(q,b\right)\in\mathcal{X}$,
	and discount factor $\beta\in\left(0,1\right)$, the value function
	$J\ubeta\left(q,b\right)$~\eqref{eq:3.1} is finite.
\end{lemma}
\begin{IEEEproof}
	For a system starting in state $(q,b)$, the queue length at time $t$ is upper-bounded as $q+M_{a}t$. Thus, the optimal cost $J\ubeta\left(q,b\right)$ is bounded as follows,
	\begin{align}
		J\ubeta\left(q,b\right) & \leq\stackrel[t=0]{\infty}{\sum}\beta^{t}\left(q+M_{a}t\right)\nonumber \\
		& =\frac{q}{1-\beta}+\frac{M_{a}\beta}{\left(1-\beta\right)^{2}}<\infty.\notag
	\end{align}
\end{IEEEproof}
Define
\begin{equation}
	\varphi_{m}\left(q\right):=\underset{b\in\left\{ p_{01},p_{11}\right\} }{\arg\min}J\ubeta\left(q,b\right),\notag
\end{equation}
where ties are broken according to some well-defined rule. Let $\left(0,\varphi_{m}\left(0\right)\right)$ be a reference state, and define
\begin{equation}
	h\ubeta\left(q,b\right):=J\ubeta\left(q,b\right)-J\ubeta\left(0,\varphi_{m}\left(0\right)\right).\label{eq:C.3}
\end{equation}

\begin{lemma}\label{lem:SennottC2}
We have $h\ubeta\left(q,b\right)\geq0$ for all $\left(q,b\right)\in\mathcal{X}$ and $\beta\in\left(0,1\right)$.
\end{lemma}
\begin{IEEEproof}
It suffices to show that the function $J\ubeta(\cdot,\cdot)$ attains its minimum value at one of the points $\left(0,p_{01}\right)$ or $\left(0,p_{11}\right)$.~From Lemma~\ref{lem:lemma4-1} we have that $J\ubeta\left(q,b\right)$ is monotonic non-decreasing in $q$, and hence its minima is attained at $(0,b)$ for some $b\in\left[p_{01},p_{11}\right]$.~From Lemma \ref{lem:lemma3} we have that $J\ubeta\left(q,b\right)$ is concave in $b$; since the minima of a concave function over a convex set is attained at an extreme point of the convex set, and the set $\{(0,p_0)$,$(0,p_1)\}$ is extreme set, this concludes the proof.
\end{IEEEproof}
\begin{lemma}\label{lem:SennottC3}
There exists a non-negative function $\eta: \cX \rightarrow \bR_+$ such that $h\ubeta\left(q,b\right) \leq \eta\left(q,b\right)$ for every $(q,b) \in \cX$ and $\beta \in (0,1)$.
\end{lemma}
\begin{IEEEproof}
	Define,
	\[
	\varphi_{M}\left(q\right)\in \underset{b\in\left\{ p_{01},p_{11}\right\} }{\arg\max}J\ubeta\left(q,b\right),
	\]
	where ties are broken according to some well-defined rule.~We have,
	\begin{align}
		& J\ubeta\left(q,b\right)  \leq J\ubeta\left(q,b;1\right)\nonumber \\
		& =q+ \kappa c(1)+\beta\stackrel[i=0]{M_{a}}{\sum}p_{i}\left(b J\ubeta\left(\cQ^+(q,i,1),p_{11}\right)\right.\nonumber \\
		& \qquad\qquad\qquad\left.\vphantom{\beta\left(b J\ubeta\left(\left(q-1\right)^{+}+i,p_{11}\right)+\right.}+\left(1-b\right)J\ubeta\left(\cQ^+(q,i,0),p_{01}\right)\right)\nonumber \\
		& \leq q+ \kappa c(1)+J\ubeta\left(q+M_{a},\varphi_{M}\left(q+M_{a}\right)\right),\label{eq:C.5}
	\end{align}
	where the first inequality follows from the definition of $J\ubeta\left(q,b\right)$,~and the second follows from Lemma \ref{lem:lemma4-1}
	and the definition of $\varphi_{M}\left(q\right)$. From (\ref{eq:C.3}), (\ref{eq:C.5}) we obtain,
	\begin{equation}
		h\ubeta\left(q,b\right)\leq q+ \kappa c(1)+h\ubeta\left(q+M_{a},\varphi_{M}\left(q+M_{a}\right)\right).\label{eq:C.6}
	\end{equation}
Since $\varphi_{M}\left(q+M_{a}\right)$ is either equal to $p_{01}$, or $p_{11}$ hence, it suffices to prove our claim for all states $\left(q,b\right)\in\widetilde{\mathcal{X}}$. This result would then easily extend to all $\left(q,b\right)\in\mathcal{X}$ from (\ref{eq:C.6}), and hence we focus exclusively on proving it for $\left(q,b\right)\in\widetilde{\mathcal{X}}$. Consider the following policy denoted by $\widetilde{\pi}$: implement $\pi_{M_{d}}$ until the state $\left(0,\varphi_{m}\left(0\right)\right)$ is reached for the first time, and thereafter use an optimal policy for the $\beta$-discounted POMDP~\eqref{def:discountedpomdp} $\pi\ust_{\beta}$ as given in Proposition~\ref{thm:SennotsTheorem}. Let $\tau$ be the first passage time from the state $\left(q,b\right)$, where $b\in\left\{ p_{01},p_{11}\right\} $,~to $\left(0,\varphi_{m}\left(0\right)\right)$ under policy $\pi_{M_{d}}$. $\tau$ is finite since from Lemma \ref{lem:irreducible} we have that $\pi_{M_{d}}$
	induces an irreducible Markov chain.~For $b\in\left\{ p_{01},p_{11}\right\} $ we have,
	
	\begin{align}
		& J\ubeta\left(q,b\right) \leq\mathbb{E}_{\widetilde{\pi}}\left[\stackrel[t=0]{\tau-1}{\sum}\beta^{t}\left(Q\left(t\right)+ \kappa u\left(t\right)\right) \bigg| X\left(0\right)=\left(q,b\right)\right]\nonumber \\
		& +\mathbb{E}_{\widetilde{\pi}}\left[\stackrel[t=\tau]{\infty}{\sum}\beta^{t}\left(Q\left(t\right)+ \kappa u\left(t\right)\right) \bigg| X\left(0\right)=\left(q,b\right)\right]\nonumber \\
		& \leq c\left(q,b\right)+E\left[\beta^{\tau}\right]J\ubeta\left(0,\varphi_{m}\left(0\right)\right)\nonumber \\
		& \leq c\left(q,b\right)+J\ubeta\left(0,\varphi_{m}\left(0\right)\right),\label{eq:C.7}
	\end{align}
	where $c\left(q,b\right)$ is the expected value of the cumulative cost incurred until the first
	passage from the state $\left(q,b\right)$ to the reference state $\left(0,\varphi_{m}\left(0\right)\right)$.
	~(\ref{eq:C.7}) is the same as,
	\begin{equation}
		h\ubeta\left(q,b\right)\leq c\left(q,b\right),\;\forall\left(q,b\right)\in\widetilde{\mathcal{X}}.\notag
	\end{equation}
	Now, since $\pi_{M_{d}}$ induces an irreducible Markov chain in $\widetilde{\mathcal{X}}$, $c\left(q,b\right)$ is finite
	for all $\left(q,b\right)\in\widetilde{\mathcal{X}}$ if and
	only if $c\left(0,\varphi_{m}\left(0\right)\right)$ is finite. Let $m\left(q,b\right)$
	be the expected time for the first passage from the state $\left(q,b\right)$
	to the reference state $\left(0,\varphi_{m}\left(0\right)\right)$ under
	the policy $\pi_{M_{d}}$.~The
	average cost incurred by $\pi_{M_{d}}$ is then equal to $\frac{c\left(0,\varphi_{m}\left(0\right)\right)}{m\left(0,\varphi_{m}\left(0\right)\right)}$~\cite[p. 75]{serfozo1981optimal}.~Since average cost of $\pi_{M_{d}}$ is finite (Lemma \ref{lem:average-cost-finite}), the quantity $\frac{c\left(0,\varphi_{m}\left(0\right)\right)}{m\left(0,\varphi_{m}\left(0\right)\right)}$, and consequently $c\left(0,\varphi_{m}\left(0\right)\right)$
	is also finite. This proves that $c\left(q,b\right)$ is finite
	for all $\left(q,b\right)\in\widetilde{\mathcal{X}}$. Hence,
	we can set $\eta\left(q,b\right)=c\left(q,b\right)$ for
	all $\left(q,b\right)\in\widetilde{\mathcal{X}}$. This concludes
	the proof.
\end{IEEEproof}
\begin{lemma}\label{lem:SennottC3star}
Let $\eta(\cdot,\cdot)$ be as given in Lemma~\ref{lem:SennottC3}.~For all $\left(q,b\right)\in\mathcal{X}$ and $j\in\left\{ 0,\ldots,M_{d}\right\} $, we have,
\begin{equation}
\underset{\left(\overline{q},\overline{b}\right)\in\mathcal{X}}{\sum}P_{\left(q,b\right),\left(\overline{q},\overline{b}\right)}\left(j\right)\cdot\eta\left(\overline{q},\overline{b}\right)<\infty,\label{eq:C.4}
	\end{equation}
	where $P_{\left(q,b\right),\left(\overline{q},\overline{b}\right)}\left(j\right)$
	is the transition probability from state $\left(q,b\right)$
	to $\left(\overline{q},\overline{b}\right)$ if action $j$ is
	taken in state $\left(q,b\right)$, and $\eta\left(q,b\right)$ is as in Lemma~\ref{lem:SennottC3}.
\end{lemma}
\begin{IEEEproof}
	The proof follows if we show that for each $\left(q,b\right)$
	and $j$, the probabilities $P_{\left(q,b\right),\left(\overline{q},\overline{b}\right)}\left(j\right)$
	are non-zero for only a finite number of $\left(\overline{q},\overline{b}\right)\in\mathcal{X}$.
	Consider the set $\mathcal{A}=\mathcal{A}_{q}\times\mathcal{A}_{b}$
	where $\times$ is the cartesian product operator, $\mathcal{A}_{q}=\left\{ \left(q-j\right)^{+},\ldots,\left(q-j\right)^{+}+M_{a}\right\} $,
	and $\mathcal{A}_{b}=\left\{ p_{01},p_{11},\mathcal{T}\left(b\right)\right\} $.
	For any given $\left(q,b\right)$ and $j$, $P_{\left(q,b\right),\left(\overline{q},\overline{b}\right)}\left(j\right)$
	is non-zero only if $\left(\overline{q},\overline{b}\right)\in\mathcal{A}$.
	Since $\mathcal{A}$ is a finite set, the LHS of (\ref{eq:C.4}) is
	finite. This concludes the proof.
\end{IEEEproof}

\end{document}